\documentclass[onecolumn]{IEEEtran}
\usepackage[latin9]{inputenc}
\usepackage{mathrsfs}
\usepackage{amsmath}
\usepackage{amsthm}
\usepackage{amssymb}
\usepackage{graphicx}
\usepackage[pdfusetitle,
 bookmarks=true,bookmarksnumbered=true,bookmarksopen=true,bookmarksopenlevel=1,
 breaklinks=false,pdfborder={0 0 0},pdfborderstyle={},backref=false,colorlinks=false]
 {hyperref}

\makeatletter
\theoremstyle{plain}
\newtheorem{notation}{\protect\notationname}
\theoremstyle{remark}
\newtheorem{rem}{\protect\remarkname}
\theoremstyle{plain}
\newtheorem{lem}{\protect\lemmaname}
\theoremstyle{plain}
\newtheorem{assumption}{\protect\assumptionname}
\theoremstyle{plain}
\newtheorem{prop}{\protect\propositionname}
\theoremstyle{plain}
\newtheorem{thm}{\protect\theoremname}
\theoremstyle{definition}
\newtheorem{defn}{\protect\definitionname}

\date{}
\usepackage{xprintlen}
\makeatother

\providecommand{\assumptionname}{Assumption}
\providecommand{\definitionname}{Definition}
\providecommand{\lemmaname}{Lemma}
\providecommand{\notationname}{Notation}
\providecommand{\propositionname}{Proposition}
\providecommand{\remarkname}{Remark}
\providecommand{\theoremname}{Theorem}

\begin{document}
\title{Distributed Newton Optimization with Maximized Convergence Rate}
\author{Damián Marelli, Yong Xu$^{\dagger}$, Minyue Fu, \emph{Fellow IEEE},
and Zenghong Huang\thanks{Damián Marelli is with the School of Automation, Guangdong University
of Technology, Guangzhou, China, and with the French Argentine International
Center for Information and Systems Sciences, National Scientific and
Technical Research Council, Argentina. Email: {\footnotesize\texttt{Damian.Marelli@newcastle.edu.au}}{\footnotesize .}}\thanks{Yong Xu is with the School of Automation, Guandong University of Technology,
China. Email:{\footnotesize{} }{\footnotesize\texttt{xuyong809@163.com.}}}\thanks{Minyue Fu is with the School of Electrical Engineering and Computer
Science, University of Newcastle, Callaghan, NSW 2308, Australia.}\thanks{Zenghong Huang is with the School of Automation, Guandong University
of Technology, China. Email:{\footnotesize{} }{\footnotesize\texttt{zenghong9527@foxmail.com.}}}\thanks{$\dagger$Corresponding author.}\thanks{This work was supported by the Argentinean Agency for Scientific and
Technological Promotion (PICT- 201-0985) and by the National Natural
Science Foundation of China (Grant Nos. 61633014, 61803101, U1911401
and U1701264).}}
\maketitle
\begin{abstract}
The distributed optimization problem is set up in a collection of
nodes interconnected via a communication network. The goal is to find
the minimizer of a global objective function formed by the addition
of partial functions locally known at each node. A number of methods
are available for addressing this problem, having different advantages.
The goal of this work is to achieve the maximum possible convergence
rate. As the first step towards this end, we propose a new method
which we show converges faster than other available options. As with
most distributed optimization methods, convergence rate depends on
a step size parameter. As the second step towards our goal we complement
the proposed method with a fully distributed method for estimating
the optimal step size that maximizes convergence speed. We provide
theoretical guarantees for the convergence of the resulting method
in a neighborhood of the solution. Also, for the case in which the
global objective function has a single local minimum, we provide a
different step size selection criterion together with theoretical
guarantees for convergence. We present numerical experiments showing
that, when using the same step size, our method converges significantly
faster than its rivals. Experiments also show that the distributed
step size estimation method achieves an asymptotic convergence rate
very close to the theoretical maximum.
\end{abstract}

\section{Introduction\protect\label{sec:Introduction}}

A networked system is a web of intelligent sensing and computing devices
connected via a communication network. Its main goal is to carry out
a computational task in a distributed manner, by executing a cooperative
strategy over all the nodes of the network without centralized coordination.
The design of distributed algorithms is constrained by the fact that
each node is limited in computational power and communication bandwidth.
Distributed algorithms are available for parameter estimation~\cite{Xiao2006,Marelli2015},
Kalman filtering~\cite{Ribeiro2010}, control~\cite{Massioni2009,DAndrea2003},
optimization~\cite{Yang2019}, etc.

The goal of a distributed optimization method is to minimize a cost
function formed by a sum of local functions which are only known by
each node~\cite{Sayed2014,Yang2019,Nedic2018,Nedich2015,Yang2010}.
It finds applications in power systems, sensor networks, smart buildings,
smart manufacturing, etc. The available distributed optimization methods
can be classified according to different criteria. We describe below
those criteria used in this work.

One classification criterion is between sequential and simultaneous
methods. In a sequential method, nodes take turns to tune its local
variables using its local function as well as information received
from other nodes~\cite{Guerbuezbalaban2015}. The main drawback of
these methods is that they do not scale well for large networks, since
many turns are needed to guarantee that all nodes are visited. Also,
for a fully distributed implementation, a distributed mechanism is
needed to guarantee that each node is regularly visited in the sequence.
A popular approach within this line are methods based on alternating
direction method of multipliers~\cite{Boyd2011,Wei2012}. In contrast
to sequential methods, a simultaneous method iterates over a computation
step, in which all nodes carry out local computations, and a communication
step, in which nodes communicate information with their neighbors.
These two steps typically depend on the number of local neighbors
and not on the network size. In this way, simultaneous methods avoid
the scalability problems of sequential ones.

Another classification criterion is between methods using first order
derivatives and those using derivatives of second order. There is
a vast literature on first order methods~\cite{nedic2009distributed,jakovetic2014fast,yuan2016convergence,shi2015extra,yang2018global,pu2020push}
as well as the survey~\cite{Yang2019}. As with centralized optimization
methods, the advantage of first order distributed methods is that
they are simpler to implement and analyze. However, second order distributed
methods converge much faster, leading to less computational and communication
requirements.

Upon convergence, a distributed optimization method needs to guarantee
that all nodes obtain the same optimal value. Another classification
criterion is based on how the distributed method guarantees this inter-node
matching property. One approach consists in adding a constraint to
the optimization program forcing this property. In~\cite{Mokhtari2016,mansoori2019fast}
the resulting constrained optimization problem is solved by adding
a penalization term. This gives an approximate solution which becomes
exact as the step size used in the optimization recursions decreases
to zero. This has the disadvantage of slowing down convergence. This
is avoided in~\cite{Tutunov2019} by solving the constrained optimization
program via its dual. However, at each iteration each node needs to
solve a local optimization problem needed to evaluate the Lagrange
dual function. These two approaches were also considered in~\cite{Jadbabaie2009}.
All these methods require the local functions at each node to be strongly
convex, which is a somehow strong requirement, since these local functions
typically represent partial information about the variables to be
optimized. They also require, at each iteration, inverting a matrix
using a recursive inversion formula. This requires running sub-iterations,
each involving a computation/communication step, between every two
main iterations. These drawbacks are avoided by forcing the inter-node
matching using average consensus. By doing so the resulting optimization
program is unconstrained. As pointed out in~\cite{Varagnolo2015},
an additional advantage of this approach is that the large literature
available for average consensus permits guaranteeing its robustness
to asynchronous communications, packet losses, time-varying network
topology, undirected communication links, etc.

In view of the above classification, we are interested in simultaneous
second order distributed methods based on average consensus. A seminal
method within this line was introduced in~\cite{Varagnolo2015,zanella2011newton}.
It was then extended in~\cite{zanella2012asynchronous} to deal with
asynchronous communications, in~\cite{bof2018multiagent} to deal
with lossy communications, and in~\cite{Li2018} to deal with event-triggered
communications in the context of continuous-time optimization methods.
A fast convergent variant of this method was proposed in~\cite{zhang2020distributed}
using finite-time average consensus. However, this kind of consensus
requires global knowledge of the network structure, hence the method
is not fully distributed. In this work we build upon this line by
focusing on convergence speed. Our contributions are the following:
(1) We propose a variant of the aforementioned distributed method.
The proposed variant converges much faster than the other methods,
for the same optimization step size. Moreover, it is also more robust
to the choice of the step size, in the sense that it converges with
step sizes which are large enough to cause the divergence of other
methods. This permits choosing larger step sizes to speed up convergence.
(2) We propose a fully distributed method for determining, at each
node, the optimal step size in the sense of maximizing the asymptotic
convergence speed. We also do a theoretical convergence analysis guaranteeing
the convergence of the resulting distributed optimization procedure,
when used in conjunction with the proposed step size selection algorithm,
in a neighborhood of the solution. (3) In the case in which the global
objective function has a single local minimum, we provide a different
step size selection criterion, together with sufficient conditions
to guarantee the global convergence of the proposed algorithm. This
result is stronger than most available ones in the sense that we require
that the global objective function, rather than each local function,
is strongly convex.

The rest of the paper is organized as follows. In Section~\ref{sec:Problem-statement}
we state the research problem. In Section~\ref{sec:proposed-method}
we derive the proposed distributed optimization method and compare
it with other available methods. In Section~\ref{sec:Step-size-local}
we describe the proposed distributed method for step size selection.
In Section~\ref{sec:Step-size-global} we do the global convergence
analysis in the case of a single local minimum. In Section~\ref{sec:Numerical-example}
we apply our results to a practical problem, namely, target localization,
and present numerical experiments confirming our claims. Concluding
remarks are given in Section~\ref{sec:Conclusion}.

\section{Problem statement\protect\label{sec:Problem-statement}}
\begin{notation}
The set of natural and real numbers are denoted by $\mathbb{N}$ and
$\mathbb{R}$, respectively. For a vector $x\in\mathbb{R}^{N}$, we
use $\left\Vert x\right\Vert $ to denote its $2$-norm and $\left\Vert x\right\Vert _{P}=\sqrt{x^{\top}Px}$
to denote its $P>0$ weighted 2-norm. For a matrix $A\in\mathbb{R}^{N\times N}$,
we use $\left\Vert A\right\Vert _{\mathrm{F}}$ to denote its Frobenius
norm, $\left\Vert A\right\Vert $ to denote its operator norm (induced
by the vector $2$-norm) and $\rho(A)$ its spectral radius. Also,
$\mathbf{I}_{N}$ denotes the $N\times N$ identity matrix and $\mathbf{1}_{N}$
the $N$-dimensional column vector filled with ones. To simplify the
notation we often omit the subscript in $\mathbf{I}_{N}$ and $\mathbf{1}_{N}$
when the dimension can be clearly inferred from the context. We use
$\mathrm{col}\left(x_{1},\cdots,x_{I}\right)$ to denote the column
vector formed by stacking the elements $x_{1},\cdots,x_{I}$ and $\otimes$
to denote the Kronecker product. We use $\leq$ to denote the non-strict
partial order on $\mathbb{R}^{N}$ defined by $x\leq y$ if $x_{i}\leq y_{i}$,
for all $i=1,\cdots,I$, and $<$ to denote the strict partial order
corresponding to $\leq$, i.e., $x<y$ if $x\leq y$ and $x\neq y$.
Finally, $O(x)$ denotes Bachmann--Landau's big O notation $O(x)$
as $x\rightarrow0$.
\end{notation}

We have a network of $I$ nodes. Node~$i$ can evaluate the function
$f^{i}:\mathbb{R}^{N}\rightarrow\mathbb{R}$ and send messages to
its out-neighbors $\mathcal{N}_{i}\subseteq\{1,\cdots,I\}$ using
a consensus network. The communication link from node~$i$ to node
$j\in\mathcal{N}_{i}$ has time-invariant gain $w^{j,i}$. We assume
that $w^{i,j}=0$ if $j\neq\mathcal{N}_{i}$. We also assume that
the graph induced by the communication network is balanced (i.e.,
possibly directed) and strongly connected. This implies that matrix
$W=\left[w^{i,j}\right]_{i,j=1}^{I}$ is doubly stochastic and primitive.
As explained in~\cite{OlfatiSaber2007}, a consequence of this is
that $\mathbf{1}_{I}^{\top}W=\mathbf{1}_{I}^{\top}$ and $W\mathbf{1}_{I}=\mathbf{1}_{I}$.
Also, for any $x_{1}=\left[x_{1}^{i},\cdots,x_{1}^{I}\right]^{\top}\in\mathbb{R}^{I}$,
the sequence generated by $x_{k+1}=Wx_{k}$ satisfies
\[
\lim_{k\rightarrow\infty}x_{k}=\mathbf{1}_{I}\otimes\frac{1}{I}\sum_{i=1}^{I}x_{1}^{i}.
\]

 The goal of distributed optimization is to design a distributed
method for solving the following minimization problem 
\begin{equation}
x_{\star}\in\underset{x\in\mathbb{R}^{N}}{\arg\min}f(x)\quad\text{with}\quad f(x)=\frac{1}{I}\sum_{i=1}^{I}f^{i}(x).\label{eq:prob}
\end{equation}
Solving~(\ref{eq:prob}) using a centralized or distributed method
requires certain assumptions on the objective function $f$, e.g.,
convexity, quasi-convexity, everywhere positive definite Hessian matrix,
etc. These assumptions may be too strong in certain applications.
When none of these assumptions can be made, it is often enough to
solve 
\begin{equation}
x_{\star}\in\underset{x\in\mathbb{R}^{N}}{\mathrm{loc\,min}}f(x),\label{eq:prob-loc}
\end{equation}
where $\mathrm{loc\,min}$ denotes the set of local minimizers of
$f$.

As mentioned in the introduction, a number of methods are available
for solving either~(\ref{eq:prob}) or~(\ref{eq:prob-loc}), and
our preferred choice is the family of methods in~\cite{Varagnolo2015,zanella2011newton,zanella2012asynchronous,bof2018multiagent,Li2018}.
In Section~\ref{sec:proposed-method} we derive a new method that
can be regarded as a variant of the aforementioned ones. As we explain,
its design aims at maximizing convergence speed. In Section~\ref{sec:Step-size-local}
we focus on the general problem~(\ref{eq:prob-loc}) in which convergence
to the global optimum cannot be guaranteed. We propose a step size
selection criterion to maximize local convergence speed. In Section~\ref{sec:Step-size-global}
we focus on~(\ref{eq:prob}) and provide a selection criterion to
guarantee convergence to the global optimum, in the case where the
global function $f$ is strongly convex.

\section{Proposed method\protect\label{sec:proposed-method}}

In this section we introduce the proposed distributed optimization
method with fixed step size. In Section~\ref{sec:Static,-dynamic-and}
we introduce some background on static and dynamic average consensus.
Using this, in Section~\ref{sec:proposed-method} we derive the proposed
method, and in Section~\ref{sec:Comparison-with-similar} describe
its differences with respect to the variant that has been used in
the literature. In Section~\ref{subsec:State-space-representation1}
we derive two state-space representations of the algorithm which are
instrumental for our convergence studies.

\subsection{Static and dynamic average consensus\protect\label{sec:Static,-dynamic-and}}

As described in Section~\ref{sec:Introduction}, we are interested
in methods achieving inter-node matching of minimization parameters
via average consensus. In this section we briefly describe the options
available for doing so.

Suppose that, in the network described in Section~\ref{sec:Problem-statement},
each node~$i$ knows a variable $v^{i}\in\mathbb{V}$, $i\in\mathbb{N}$,
where $\mathbb{V}$ is a vector space. In order to make the presentation
valid in the general case, we assume that $\mathbb{V}$ is an arbitrary
vector space, i.e., each $v^{i}$ can be either a scalar, a vector,
a matrix, etc. The goal of (static) average consensus is to compute
the average $u=\frac{1}{I}\sum_{i=1}^{I}v^{i}$ in a distributed manner.
This is done using the following iterations
\[
v_{k+1}^{i}=\sum_{j=1}^{I}w^{i,j}v_{k}^{j},
\]
initialized by $v_{1}^{i}=v^{i}$~\cite{OlfatiSaber2007}. Letting
$\mathbf{v}_{k}=\mathrm{col}\left(v_{k}^{1},\cdots,v_{k}^{I}\right)\in\mathbb{V}^{I}$
we can write the above compactly as follows
\begin{equation}
\mathbf{v}_{k+1}=W\mathbf{v}_{k}.\label{eq:FOC}
\end{equation}

Suppose now that each node~$i$ knows a time-varying sequence of
variables $v_{k}^{i}\in\mathbb{V}$, $k\in\mathbb{N}$. The goal of
the dynamic average consensus technique~\cite{zhu2010discrete} is
to obtain, at each $k$, an estimate of the average $u_{k}=\frac{1}{I}\sum_{i=1}^{I}v_{k}^{i}$.
This is done as follows: Suppose that at time $k$ node~$i$ knows
an estimate $u_{k-1}^{i}$ of $u_{k-1}$. It then transmits the following
\emph{message}
\begin{equation}
s_{k}^{i}=u_{k-1}^{i}+v_{k}^{i}-v_{k-1}^{i},\label{eq:tx}
\end{equation}
to its out-neighbors. On reception, node $i$ obtains
\begin{equation}
u_{k}^{i}=\sum_{j=1}^{I}w^{i,j}s_{k}^{j}.\label{eq:rx}
\end{equation}
The above iterations are initialized by $s_{1}^{i}=v_{1}^{i}$. We
can combine~(\ref{eq:tx})-(\ref{eq:rx}) in two ways, namely, in
\emph{message form}
\[
s_{k+1}^{i}=\sum_{j=1}^{I}w^{i,j}s_{k}^{j}+v_{k+1}^{i}-v_{k}^{i},
\]
or in \emph{estimate form}
\[
u_{k+1}^{i}=\sum_{j=1}^{I}w^{i,j}\left(u_{k}^{j}+v_{k+1}^{j}-v_{k}^{j}\right).
\]

\subsection{The proposed method}

 The essential idea consists in distributing the Newton iterations
\begin{equation}
x_{k+1}=x_{k}-\alpha_{k}\left[\nabla^{2}f\left(x_{k}\right)\right]^{-1}\nabla f\left(x_{k}\right),\label{eq:newton}
\end{equation}
where $\alpha_{k}$ is called the step size at time $k$. To this
end we make use of the dynamic average consensus technique~\cite{zhu2010discrete}.
Let 
\[
\bar{f}\left(x_{k}^{1},\cdots,x_{k}^{I}\right)=\frac{1}{I}\sum_{i=1}^{I}f^{i}\left(x_{k}^{i}\right).
\]
We can obtain an estimate $g_{k}^{i}$ of $\nabla\bar{f}\left(x_{k}^{1},\cdots,x_{k}^{I}\right)=\frac{1}{I}\sum_{i=1}^{I}\nabla f^{i}\left(x_{k}^{i}\right)$,
at each node~$i$, by applying dynamic average consensus on the inputs
$\nabla f^{i}\left(x_{k}^{i}\right)$. This yields the following recursions
written in estimate form
\[
g_{k}^{i}=\sum_{j=1}^{I}w^{i,j}\left[g_{k-1}^{j}+\nabla f^{j}\left(x_{k}^{j}\right)-\nabla f^{j}\left(x_{k-1}^{j}\right)\right].
\]
We can do the same to obtain an estimate $H_{k}^{i}$ of the Hessian
$\nabla^{2}\bar{f}\left(x_{k}^{1},\cdots,x_{k}^{I}\right)=\frac{1}{I}\sum_{i=1}^{I}\nabla^{2}f^{i}\left(x_{k}^{i}\right)$
using the inputs $\nabla^{2}f^{i}\left(x_{k}^{i}\right)$. This gives
\[
H_{k}^{i}=\sum_{j=1}^{I}w^{i,j}\left[H_{k-1}^{j}+\nabla^{2}f^{j}\left(x_{k}^{j}\right)-\nabla^{2}f^{j}\left(x_{k-1}^{j}\right)\right].
\]

It is easy to see that $\bar{f}\left(x,\cdots,x\right)=f(x)$. Hence,
if 
\begin{equation}
x_{k}^{i}\simeq x_{k}\quad\text{for all}\quad i\in\{1,\cdots,I\},\label{eq:match}
\end{equation}
and some $x_{k}$, then $g_{k}^{i}$ and $H_{k}^{i}$ are estimates
of $\nabla f\left(x_{k}\right)$ and $\nabla^{2}f\left(x_{k}\right)$,
respectively. Thus, in principle, each node~$i$ could use $g_{k}^{i}$
and $H_{k}^{i}$, in place of $\nabla f\left(x_{k}\right)$ and $\nabla^{2}f\left(x_{k}\right)$
to locally carry out the iterations~(\ref{eq:newton}). This would
yield, at node~$i$, the following sequence of estimates of $x_{\star}$
\[
\breve{x}_{k+1}^{i}=\breve{x}_{k}^{i}-\alpha_{k}^{i}\left[H_{k}^{i}\right]^{-1}g_{k}^{i}.
\]
But the above requires~(\ref{eq:match}) to hold, with $x_{k}^{i}$
replaced by $\breve{x}_{k}^{i}$. In order to enforce that, once again
we apply dynamic average consensus on the inputs $\breve{x}_{k}^{i}$.
Writing the result in message form we obtain
\[
x_{k+1}=\sum_{j=1}^{I}w^{i,j}x_{k}^{j}+\breve{x}_{k+1}^{i}-\breve{x}_{k}^{i}=\sum_{j=1}^{I}w^{i,j}x_{k}^{j}-\alpha_{k}^{i}\left[H_{k}^{i}\right]^{-1}g_{k}^{i}.
\]

Finally, since the approximations $H_{k}^{i}$ to the Hessian are
obtained via dynamic average consensus on the local Hessian matrices
$\nabla f^{j}\left(x_{k}^{j}\right)$, and the latter may fail to
be positive definite, some mechanism is required to guarantee that
$H_{k}^{i}$ is positive definite. To do so we let $\beta>0$ and
use a map $B:$ to guarantee that $B\left(H_{k}^{i}\right)\geq\beta^{-1}\mathbf{I}$.
The map $B$ is defined as follows: Let $H\in\mathbb{R}^{N}$ be symmetric
and $H=U\Lambda U^{\top}$ be its spectral decomposition, with $\Lambda=\mathrm{diag}\left(\lambda_{1},\cdots,\lambda_{N}\right)$.
Let $\tilde{\Lambda}=\mathrm{diag}\left(\tilde{\lambda}_{1},\cdots,\tilde{\lambda}_{N}\right)$
with $\tilde{\lambda}_{i}=\lambda_{i}$ if $\lambda_{i}\geq\beta^{-1}$
and $\tilde{\lambda}_{i}=\beta^{-1}$ otherwise. Then 
\[
B(H)=U\tilde{\Lambda}U^{\top}.
\]
The question then arises as to how to choose the parameter $\beta$.
This is given in Assumptions~\ref{assu:local} and~\ref{assu:global}
of our main results given in Sections~\ref{sec:Step-size-local}
and~\ref{sec:Step-size-global}, respectively.

To summarize the above, the proposed algorithm is given by the following
recursions, 
\begin{align}
x_{k+1}^{i} & =\sum_{j=1}^{I}w^{i,j}x_{k}^{j}-\alpha_{k}^{i}B\left(H_{k}^{i}\right)^{-1}g_{k}^{i},\label{eq:dist-method-1}\\
g_{k+1}^{i} & =\sum_{j=1}^{I}w^{i,j}\left[g_{k}^{j}+\nabla f^{j}\left(x_{k+1}^{j}\right)-\nabla f^{j}\left(x_{k}^{j}\right)\right],\label{eq:dist-method-2}\\
H_{k+1}^{i} & =\sum_{j=1}^{I}w^{i,j}\left[H_{k}^{j}+\nabla^{2}f^{j}\left(x_{k+1}^{j}\right)-\nabla^{2}f^{j}\left(x_{k}^{j}\right)\right].\label{eq:dist-method-3}
\end{align}
which are initialized by
\begin{equation}
x_{1}^{i}=x_{\mathrm{init}}^{i},\quad g_{1}^{i}=\nabla f^{i}\left(x_{\mathrm{init}}^{i}\right)\quad\text{and}\quad H_{1}^{i}=\nabla^{2}f^{i}\left(x_{\mathrm{init}}^{i}\right).\label{eq:init-cond}
\end{equation}

\begin{rem}
The reader may wonder why we choose to write~(\ref{eq:dist-method-1})
in message form while~(\ref{eq:dist-method-2}) and~(\ref{eq:dist-method-3})
in estimate form. This is done to put the algorithm in a form compatible
with other algorithms in the literature.
\end{rem}
\begin{rem}
Notice that the information exchanged by each node at each time step
does not grow with the network size, as it depends on the number of
out neighbors of each node. This property is common to all simultaneous,
second order methods based on average consensus~\cite{Varagnolo2015,zanella2011newton,zanella2012asynchronous,bof2018multiagent,Li2018}.
\end{rem}

\subsection{Comparison with similar algorithms\protect\label{sec:Comparison-with-similar}}

As mentioned above, the proposed algorithm~(\ref{eq:dist-method-1})-(\ref{eq:dist-method-3})
is a variant of the algorithm used in~\cite{Varagnolo2015,zanella2011newton,zanella2012asynchronous,bof2018multiagent,Li2018}.
The latter differ from~(\ref{eq:dist-method-1})-(\ref{eq:dist-method-3})
in essentially two aspects. The first one is that consensus is not
done on the parameters $x_{k}^{i}$. This means that~(\ref{eq:dist-method-1})
is replaced by
\begin{align}
x_{k+1}^{i} & =x_{k}^{i}-\alpha_{k}^{i}B\left(H_{k}^{i}\right)^{-1}g_{k}^{i}.\label{eq:algA}
\end{align}
Together with~(\ref{eq:dist-method-2})-(\ref{eq:dist-method-3}),
equation~(\ref{eq:algA}) forms an algorithm that, for latter reference,
we refer to as Algorithm~A.

The second difference consists in using the following transformation
of~(\ref{eq:newton}) 
\begin{equation}
x_{k+1}=(1-\alpha)x_{k}+\alpha_{k}^{i}\left(\nabla^{2}f\left(x_{k}\right)\right)^{-1}\ell\left(x_{k}\right),\label{eq:trick}
\end{equation}
where $\ell(x)=\nabla^{2}f(x)x-\nabla f(x)$. Using dynamic average
consensus, we can estimate $\ell\left(x_{k}\right)$ at each node
using
\begin{equation}
l_{k+1}^{i}=\sum_{j=1}^{I}w^{i,j}\left[l_{k}^{j}+\ell^{j}\left(x_{k+1}^{j}\right)-\ell^{j}\left(x_{k}^{j}\right)\right],\label{eq:algB2}
\end{equation}
where $\ell^{i}(x)=\nabla^{2}f^{i}(x)x-\nabla f^{i}(x)$. We can then
distribute iterations~(\ref{eq:trick}) as follows
\begin{align}
x_{k+1}^{i} & =\left(1-\alpha_{k}^{i}\right)\sum_{j=1}^{I}w^{i,j}x_{k}^{j}+\alpha_{k}^{i}B\left(H_{k}^{i}\right)^{-1}l_{k}^{i}.\label{eq:algB1}
\end{align}
We refer to the algorithm resulting from~(\ref{eq:algB1}),~(\ref{eq:algB2})
and~(\ref{eq:dist-method-3}) as Algorithm~B.

Finally, the algorithm used in~\cite{Varagnolo2015,zanella2011newton,zanella2012asynchronous,bof2018multiagent,Li2018},
apart from other minor differences, essentially consists in combining
the modifications introduced by Algorithms~A and~B. This leads to
the recursions formed by
\begin{equation}
x_{k+1}^{i}=\left(1-\alpha_{k}^{i}\right)x_{k}^{i}+\alpha_{k}^{i}B\left(H_{k}^{i}\right)^{-1}l_{k}^{i},\label{eq:VZCPS}
\end{equation}
together with~(\ref{eq:dist-method-3}) and~(\ref{eq:algB2}). We
refer to it as Algorithm~VZCPS, standing for the initials of the
authors which proposed it.

As we show with experiments in Section~\ref{sec:Numerical-example},
the modifications~(\ref{eq:algA}) and~(\ref{eq:algB2})-(\ref{eq:algB1}),
introduced by Algorithms~A and~B, respectively, drastically slow
down convergence and can cause instability. More precisely,~(\ref{eq:algA})
does not guarantee the convergence of each $x_{k}^{i}$ to $x_{\star}$,
due to the lack of consensus on these parameters. Convergence only
occurs if modification~(\ref{eq:algB2})-(\ref{eq:algB1}) is also
considered, i.e., in the VZCPS algorithm, although at a much slower
rate. However, a feature of the latter modification is that the first
term in~(\ref{eq:algB1}) pushes the local variables $x_{k}^{i}$
towards zero at each iteration. This pushing is compensated by the
second term, but only after consensus on the parameters $x_{k}^{i}$
is reached. Before this happens, this zero pushing effect has a negative
influence if the minimizing parameters $x_{\star}$ are far from zero.
As we show in Section~\ref{sec:Numerical-example}, this can slow
down convergence and even cause instability.

\subsection{State-space representation\protect\label{subsec:State-space-representation1}}

We introduce the following required notation.
\begin{notation}
Let $A=\frac{1}{I}\mathbf{1}_{I}\mathbf{1}_{I}^{\top}$, $\mathbf{A}\mathbf{=}A\otimes\mathbf{I}_{N}$,
$\mathbf{a}=\frac{1}{I}\mathbf{1}_{I}\otimes\mathbf{I}_{N}$ and $\tilde{\mathbf{I}}\mathbf{=}\mathbf{I}-\mathbf{A}$.
Let also $\mathbf{x}_{\star}=\mathbf{1}_{I}\otimes x_{\star}$, $\mathbf{x}_{k}=\mathrm{col}\left(x_{k}^{1},\cdots,x_{k}^{I}\right)$,
$\bar{\mathbf{x}}_{k}=\mathbf{A}\mathbf{x}_{k}$ and $\tilde{\mathbf{x}}_{k}=\mathbf{x}_{k}-\bar{\mathbf{x}}_{k}$.
We similarly define $\mathbf{g}_{k}=\mathrm{col}\left(g_{k}^{1},\cdots,g_{k}^{I}\right)$,
$\bar{\mathbf{g}}_{k}=\mathbf{A}\mathbf{g}_{k}$ and $\tilde{\mathbf{g}}_{k}=\mathbf{g}_{k}-\bar{\mathbf{g}}_{k}$
as well as $\mathbf{h}_{k}=\mathrm{col}\left(H_{k}^{1},\cdots,H_{k}^{I}\right)$,
$\bar{\mathbf{h}}_{k}=\mathbf{A}\mathbf{h}_{k}$, $\tilde{\mathbf{h}}_{k}=\mathbf{h}_{k}-\bar{\mathbf{h}}_{k}$,
$\mathbf{h}_{\star}=\mathbf{1}_{I}\otimes\nabla^{2}f\left(x_{\star}\right)$
and $\mathbf{H}_{\star}=\mathbf{I}_{I}\otimes\nabla^{2}f\left(x_{\star}\right)$.
Finally, $\boldsymbol{\alpha}_{k}=\mathrm{diag}\left(\alpha_{k}^{1},\cdots,\alpha_{k}^{I}\right)\otimes\mathbf{I}_{I}$.
\end{notation}
\begin{rem}
In the above notation, $\bar{\mathbf{x}}_{k}$ is a block vector with
all its sub-vectors equal to the average 
\begin{equation}
\bar{x}_{k}=\frac{1}{I}\sum_{i=1}^{I}x_{k}^{i}.\label{eq:x-ave}
\end{equation}
Also, $\tilde{\mathbf{x}}_{k}$ is a block vector whose $i$-th sub-vector
is given by $\tilde{x}_{k}^{i}=x_{k}^{i}-\bar{x}_{k}$. Finally, notice
that $\mathbf{h}_{k}$ is a (column) vector of matrices, i.e., $\mathbf{h}_{k}\in\mathbb{R}^{IN\times N}$.
\end{rem}
\begin{notation}
\label{nota:1}For a block vector $\mathbf{x}=\mathrm{col}\left(x^{1},\cdots,x^{I}\right),$let
\begin{align*}
\mathfrak{g}\left(\mathbf{x}\right) & =\mathrm{col}\left(\nabla f^{1}\left(x^{1}\right),\cdots,\nabla f^{I}\left(x^{I}\right)\right),\\
\mathfrak{h}\left(\mathbf{x}\right) & =\mathrm{col}\left(\nabla^{2}f^{1}\left(x^{1}\right),\cdots,\nabla^{2}f^{I}\left(x^{I}\right)\right),\\
\mathfrak{H}\left(\mathbf{x}\right) & =\mathrm{diag}\left(\mathfrak{h}\left(\mathbf{x}\right)\right),
\end{align*}
and for a block diagonal matrix $\mathbf{H}=\mathrm{diag}\left(H^{1},\dots.H^{I}\right)$,
let
\[
\mathfrak{B}\left(\mathbf{H}\right)=\mathrm{diag}\left(B\left(H^{1}\right),\cdots,B\left(H^{I}\right)\right).
\]
Let also $\mathbf{H}_{k}=\mathrm{diag}\left(\mathbf{h}_{k}\right)$
and $\mathbf{B}_{k}=\mathfrak{B}\left(\mathbf{H}_{k}\right)$. Let
finally $W=\left[w^{i,j}\right]_{i,j=1}^{I}$ and $\mathbf{W}=W\otimes\mathbf{I}_{N}$.
\end{notation}
Using the above notation we can write~(\ref{eq:dist-method-1})-(\ref{eq:dist-method-3})
in the following block state-space form
\begin{align}
\mathbf{x}_{k+1} & =\mathbf{W}\mathbf{x}_{k}-\boldsymbol{\alpha}_{k}\mathbf{B}_{k}^{-1}\mathbf{g}_{k},\label{eq:ssa1}\\
\mathbf{g}_{k+1} & =\mathbf{W}\left[\mathbf{g}_{k}+\mathfrak{g}\left(\mathbf{x}_{k+1}\right)-\mathfrak{g}\left(\mathbf{x}_{k}\right)\right],\label{eq:ssa2}\\
\mathbf{h}_{k+1} & =\mathbf{W}\left[\mathbf{h}_{k}+\mathfrak{h}\left(\mathbf{x}_{k+1}\right)-\mathfrak{h}\left(\mathbf{x}_{k}\right)\right].\label{eq:ssa3}
\end{align}
A problem of the above model for studying stability is that $\rho\left(\mathbf{W}\right)=1$.
Our next step is to transform~(\ref{eq:ssa1})-(\ref{eq:ssa3}) into
an equivalent model which avoids this drawback.
\begin{notation}
Let $\bar{\mathfrak{g}}\left(\mathbf{x}\right)=\mathbf{A}\mathfrak{g}\left(\mathbf{x}\right)$
and $\bar{\mathfrak{h}}\left(\mathbf{x}\right)=\mathbf{A}\mathfrak{h}\left(\mathbf{x}\right)$.
Let also $\bar{\mathfrak{H}}\left(\mathbf{x}\right)=\mathrm{diag}\left(\bar{\mathfrak{h}}\left(\mathbf{x}\right)\right)$,
$\bar{\mathbf{H}}_{k}=\mathrm{diag}\left(\bar{\mathbf{h}}_{k}\right)$,
$\tilde{\mathbf{H}}_{k}=\mathrm{diag}\left(\tilde{\mathbf{h}}_{k}\right)$
and $\tilde{\mathbf{W}}=\mathbf{W}-\mathbf{A}$.
\end{notation}
\begin{lem}
\label{lem:g-h-bar}The following equivalences hold
\begin{equation}
\bar{\mathbf{g}}_{k}=\bar{\mathfrak{g}}\left(\mathbf{x}_{k}\right),\qquad\bar{\mathbf{h}}_{k}=\bar{\mathfrak{h}}\left(\mathbf{x}_{k}\right),\qquad\bar{\mathbf{H}}_{k}=\bar{\mathfrak{H}}\left(\mathbf{x}_{k}\right).\label{eq:ghH-bar}
\end{equation}
\end{lem}
\begin{IEEEproof}
Since matrix $W$ is doubly stochastic, it follows that $\mathbf{A}\mathbf{W}=\mathbf{A}$.
Then, from~(\ref{eq:ssa2})
\begin{equation}
\bar{\mathbf{g}}_{k+1}=\mathbf{A}\mathbf{g}_{k+1}=\bar{\mathbf{g}}_{k}+\bar{\mathfrak{g}}\left(\mathbf{x}_{k+1}\right)-\bar{\mathfrak{g}}\left(\mathbf{x}_{k}\right).\label{eq:pr1}
\end{equation}
The first equation in~(\ref{eq:ghH-bar}) then follows from~(\ref{eq:pr1})
since $\bar{\mathbf{g}}_{1}=\mathbf{A}\mathbf{g}_{1}=\mathbf{A}\mathfrak{g}\left(\mathbf{x}_{1}\right)=\bar{\mathfrak{g}}\left(\mathbf{x}_{1}\right)$.
The other two equations follow using the same argument.
\end{IEEEproof}
Using Lemma~\ref{lem:g-h-bar} we can write model~(\ref{eq:ssa1})-(\ref{eq:ssa3})
as follows
\begin{align}
\bar{x}_{k+1} & =\bar{x}_{k}-\mathbf{a}^{\top}\boldsymbol{\alpha}_{k}\mathbf{B}_{k}^{-1}\mathbf{g}_{k},\label{eq:bsys1}\\
\tilde{\mathbf{x}}_{k+1} & =\tilde{\mathbf{W}}\tilde{\mathbf{x}}_{k}-\tilde{\mathbf{I}}\boldsymbol{\alpha}_{k}\mathbf{B}_{k}^{-1}\mathbf{g}_{k},\\
\tilde{\mathbf{g}}_{k+1} & =\tilde{\mathbf{W}}\left[\tilde{\mathbf{g}}_{k}+\mathfrak{g}\left(\mathbf{x}_{k+1}\right)-\mathfrak{g}\left(\mathbf{x}_{k}\right)\right],\label{eq:bsys3}\\
\tilde{\mathbf{h}}_{k+1} & =\tilde{\mathbf{W}}\left[\tilde{\mathbf{h}}_{k}+\mathfrak{h}\left(\mathbf{x}_{k+1}\right)-\mathfrak{h}\left(\mathbf{x}_{k}\right)\right],\label{eq:bsys4}
\end{align}
where $\mathbf{g}_{k}=\bar{\mathfrak{g}}\left(\mathbf{x}_{k}\right)+\tilde{\mathbf{g}}_{k}$
and $\mathbf{B}_{k}=\mathfrak{B}\left(\bar{\mathfrak{H}}\left(\mathbf{x}_{k}\right)+\tilde{\mathbf{H}}_{k}\right)$,
$\mathbf{x}_{k}=\bar{\mathbf{x}}_{k}+\tilde{\mathbf{x}}_{k}$ and
$\bar{\mathbf{x}}_{k}=\mathbf{1}_{I}\otimes\bar{x}_{k}$.

\section{Step size selection for fast local convergence\protect\label{sec:Step-size-local}}

In this section we consider the problem~(\ref{eq:prob-loc}), in
which the objective function may have a number of local minima. In
Section~\ref{subsec:non-adaptive} we derive a criterion for choosing
the step size at each node to maximize the convergence speed in a
neighborhood of the solution. This method is based on certain approximation.
In order to obtain a more accurate method, in Section~\ref{subsec:Adaptive}
we propose a distributed adaptive method estimate the optimal step
size at each node, and analyze the local stability of the resulting
distributed optimization algorithm with adaptive step size selection.

\subsection{Offline step size selection\protect\label{subsec:non-adaptive}}

We introduce the following notation.
\begin{notation}
Let $\breve{\mathbf{x}}_{k}=\mathbf{x}_{k}-\mathbf{x}_{\star}$, $\breve{\mathbf{h}}_{k}=\mathbf{h}_{k}-\mathbf{h}_{\star}$,
$\boldsymbol{\xi}_{k}=\left(\breve{\mathbf{x}}_{k},\mathbf{g}_{k},\breve{\mathbf{h}}_{k}\right)$
and $\left\Vert \boldsymbol{\xi}_{k}\right\Vert ^{2}=\left\Vert \breve{\mathbf{x}}_{k}\right\Vert ^{2}+\left\Vert \mathbf{g}_{k}\right\Vert ^{2}+\left\Vert \breve{\mathbf{h}}_{k}\right\Vert _{\mathrm{F}}^{2}$.
Let also $\boldsymbol{\Phi},\boldsymbol{\Psi}:\mathbb{M}\rightarrow\mathbb{M}$
be the linear maps with matrix representation
\[
\boldsymbol{\Phi}=\left[\begin{array}{ccc}
\mathbf{W} & 0 & 0\\
\mathbf{W}\mathfrak{H}\left(\mathbf{x}_{\star}\right)\left(\mathbf{W}-\mathbf{I}\right) & \mathbf{W} & 0\\
\mathbf{W}\mathscr{D}\mathfrak{H}\left(\mathbf{x}_{\star}\right)\left(\left(\mathbf{W}-\mathbf{I}\right)\cdot\right) & 0 & \mathbf{W}
\end{array}\right],\qquad\boldsymbol{\Psi}=\left[\begin{array}{ccc}
0 & -\mathbf{H}_{\star}^{-1} & 0\\
0 & -\mathbf{W}\mathfrak{H}\left(\mathbf{x}_{\star}\right)\mathbf{H}_{\star}^{-1} & 0\\
0 & -\mathbf{W}\mathscr{D}\mathfrak{H}\left(\mathbf{x}_{\star}\right)\left(\mathbf{H}_{\star}^{-1}\cdot\right) & 0
\end{array}\right],
\]
where $\mathscr{D}\mathfrak{h}\left(\mathbf{x}_{\star}\right)\left(\mathbf{M}\cdot\right)$
denotes the linear operator $\mathbf{y}\mapsto\mathscr{D}\mathfrak{h}\left(\mathbf{x}_{\star}\right)\left(\mathbf{M}\mathbf{y}\right)$,
with $\mathscr{\mathscr{D}\mathfrak{h}\left(\mathbf{x}_{\star}\right)}$
denoting the Fréchet derivative of $\mathfrak{h}$ at $\mathbf{x}_{\star}$.
\end{notation}
Our first result is given in Proposition~\ref{prop:localdyn1}. It
states the local linear dynamics in a neighborhood of the local optimum
$x_{\star}$.
\begin{assumption}
\label{assu:local}$\left\Vert \left[\nabla^{2}f\left(x_{\star}\right)\right]^{-1}\right\Vert \leq\beta$.
\end{assumption}
\begin{prop}
\label{prop:localdyn1} Under Assumption~\ref{assu:local}, if $\alpha_{k}^{i}=\alpha$,
for all $k\in\mathbb{N}$ and $i=1,\cdots,I$, then 
\begin{equation}
\boldsymbol{\xi}_{k+1}=\left(\boldsymbol{\Phi}+\alpha\boldsymbol{\Psi}\right)\boldsymbol{\xi}_{k}+O\left(\left\Vert \boldsymbol{\xi}_{k}\right\Vert ^{2}\right).\label{eq:loc-dyn}
\end{equation}
\end{prop}
\begin{IEEEproof}
With some abuse of notation we use $\boldsymbol{\xi}_{k+1}\left(\boldsymbol{\xi}_{k}\right)$
to denote the map $\boldsymbol{\xi}_{k}\mapsto\boldsymbol{\xi}_{k+1}$
induced by~(\ref{eq:ssa1})-(\ref{eq:ssa3}). Doing a Taylor expansion
of this map, using Fréchet derivatives, around $\boldsymbol{\xi}_{\star}=\left(0,0,0\right)$,
we obtain 
\[
\boldsymbol{\xi}_{k+1}\left(\boldsymbol{\xi}_{k}\right)=\boldsymbol{\xi}_{k+1}\left(\boldsymbol{\xi}_{\star}\right)+\mathscr{D}\boldsymbol{\xi}_{k+1}\left(\boldsymbol{\xi}_{\star}\right)\left(\boldsymbol{\xi}_{k}\right)+O\left(\left\Vert \boldsymbol{\xi}_{k}\right\Vert ^{2}\right).
\]
Clearly, $\boldsymbol{\xi}_{k+1}\left(\boldsymbol{\xi}_{\star}\right)=0$.
Also 
\[
\mathscr{D}\boldsymbol{\xi}_{k+1}\left(\boldsymbol{\xi}_{\star}\right)\left(\boldsymbol{\xi}_{k}\right)=\left(\mathscr{D}\breve{\mathbf{x}}_{k+1}\left(\boldsymbol{\xi}_{\star}\right)\left(\boldsymbol{\xi}_{k}\right),\mathscr{D}\mathbf{g}_{k+1}\left(\boldsymbol{\xi}_{\star}\right)\left(\boldsymbol{\xi}_{k}\right),\mathscr{D}\breve{\mathbf{h}}_{k+1}\left(\boldsymbol{\xi}_{\star}\right)\left(\boldsymbol{\xi}_{k}\right)\right).
\]
The result then follows since
\begin{align*}
\mathscr{D}\breve{\mathbf{x}}_{k+1}\left(\boldsymbol{\xi}_{\star}\right)\left(\boldsymbol{\xi}_{k}\right) & =\mathbf{W}\breve{\mathbf{x}}_{k}-\boldsymbol{\alpha}\mathbf{H}_{\star}^{-1}\mathbf{g}_{k},
\end{align*}
and
\begin{align*}
\mathscr{D}\mathbf{g}_{k+1}\left(\boldsymbol{\xi}_{\star}\right)\left(\boldsymbol{\xi}_{k}\right) & =\mathbf{W}\left[\mathbf{g}_{k}+\mathfrak{H}\left(\mathbf{x}_{\star}\right)\left(\mathscr{D}\breve{\mathbf{x}}_{k+1}\left(\boldsymbol{\xi}_{\star}\right)\left(\boldsymbol{\xi}_{k}\right)-\breve{\mathbf{x}}_{k}\right)\right]\\
 & =\mathbf{W}\mathfrak{H}\left(\mathbf{x}_{\star}\right)\left(\mathbf{W}-\mathbf{I}\right)\breve{\mathbf{x}}_{k}+\mathbf{W}\left(\mathbf{I}-\mathfrak{H}\left(\mathbf{x}_{\star}\right)\boldsymbol{\alpha}\mathbf{H}_{\star}^{-1}\right)\mathbf{g}_{k}
\end{align*}
and
\[
\mathscr{D}\breve{\mathbf{h}}_{k+1}\left(\boldsymbol{\xi}_{\star}\right)\left(\boldsymbol{\xi}_{k}\right)=\mathbf{W}\mathscr{D}\mathfrak{H}\left(\mathbf{x}_{\star}\right)\left(\left(\mathbf{W}-\mathbf{I}\right)\breve{\mathbf{x}}_{k}\right)-\mathbf{W}\mathscr{D}\mathfrak{H}\left(\mathbf{x}_{\star}\right)\left(\boldsymbol{\alpha}\mathbf{H}_{\star}^{-1}\mathbf{g}_{k}\right)+\mathbf{W}\breve{\mathbf{h}}_{k}.
\]
\end{IEEEproof}
The first order local dynamics~(\ref{prop:localdyn1}) reveals that
$\breve{\mathbf{h}}_{k}$ does not act as an input of neither $\breve{\mathbf{x}}_{k}$
nor $\mathbf{g}_{k}$. Hence, to study the local convergence speed
we can focus on the local dynamics of the pair $\left(\breve{\mathbf{x}}_{k},\mathbf{g}_{k}\right)$.
This is determined by the following matrix 
\[
\boldsymbol{\Gamma}(\alpha)=\left[\begin{array}{cc}
\mathbf{W} & 0\\
\mathbf{W}\mathfrak{H}\left(\mathbf{x}_{\star}\right)\left(\mathbf{W}-\mathbf{I}\right) & \mathbf{W}
\end{array}\right]+\alpha\left[\begin{array}{cc}
0 & -\mathbf{H}_{\star}^{-1}\\
0 & -\mathbf{W}\mathfrak{H}\left(\mathbf{x}_{\star}\right)\mathbf{H}_{\star}^{-1}
\end{array}\right].
\]
Clearly, when $\alpha=0$, the spectrum of $\boldsymbol{\Gamma}(0)$
consists of the eigenvalues of $W$, each having multiplicity $2N$.
Our main result describes the behavior of each of these eigenvalues
when $\alpha$ is small.
\begin{thm}
\label{thm:eigenapprox} Let $\mu_{0}\in\sigma$$(W)$ ($\sigma(X)$
denotes the spectrum of matrix $X$) and $u$ and $v$ be, respectively,
a right and left eigenvector of $W$, associated with $\mu_{0}$.
Under Assumption~\ref{assu:local}, for $\alpha>0$, there exists
$\mu\in\sigma\left(\Gamma(\alpha)\right)$ and $y\in\mathbb{R}^{N}$,
with $\left\Vert y\right\Vert =1$, satisfying 
\begin{equation}
\mu^{2}-\mu_{0}\left(2-\alpha s\right)\mu+\mu_{0}\left(\mu_{0}-\alpha s\right)\simeq0,\label{eq:sol}
\end{equation}
where $s=y^{\top}Ry$, with 
\begin{equation}
R=\frac{1}{v^{\top}u}\sum_{i=1}^{I}v^{i}u^{i}\nabla^{2}f^{i}\left(x_{\star}\right)\left[\nabla^{2}f\left(x_{\star}\right)\right]^{-1}.\label{eq:R}
\end{equation}
\end{thm}
\begin{IEEEproof}
Let and $\boldsymbol{\Gamma}(\alpha)\left[\breve{\mathbf{x}}^{\top},\mathbf{g}^{\top}\right]^{\top}=\mu\left[\breve{\mathbf{x}}^{\top},\mathbf{g}^{\top}\right]^{\top}$.
We then have
\begin{align}
\mu\breve{\mathbf{x}} & =\mathbf{W}\breve{\mathbf{x}}-\alpha\mathbf{H}_{\star}^{-1}\mathbf{g},\label{eq:gev1}\\
\mu\mathbf{g} & =\mathbf{W}\mathfrak{H}\left(\mathbf{x}_{\star}\right)\left(\mathbf{W}-\mathbf{I}\right)\breve{\mathbf{x}}+\left(\mathbf{W}-\alpha\mathbf{W}\mathfrak{H}\left(\mathbf{x}_{\star}\right)\mathbf{H}_{\star}^{-1}\right)\mathbf{g}.\label{eq:gev2}
\end{align}
Then, since $\mathbf{H}_{\star}$ and $\mathbf{W}$ commute, from~(\ref{eq:gev1})
we obtain
\[
-\alpha^{-1}\mathbf{H}_{\star}\left(\mu\mathbf{I}-\mathbf{W}\right)\breve{\mathbf{x}}=\mathbf{g}
\]
Let
\[
\mathbf{y}=-\alpha^{-1}\mathbf{H}_{\star}\breve{\mathbf{x}}.
\]
Then, 
\[
\left(\mu\mathbf{I}-\mathbf{W}\right)\mathbf{y}=\mathbf{g}
\]
and from~(\ref{eq:gev2}),
\[
\left[\mu\mathbf{I}-\left(\mathbf{W}-\alpha\mathbf{W}\mathfrak{H}\left(\mathbf{x}_{\star}\right)\mathbf{H}_{\star}^{-1}\right)\right]\left(\mu\mathbf{I}-\mathbf{W}\right)\mathbf{y}=-\alpha\mathbf{W}\mathfrak{H}\left(\mathbf{x}_{\star}\right)\mathbf{H}_{\star}^{-1}\left(\mathbf{W}-\mathbf{I}\right)\mathbf{y}.
\]
Letting $\mathbf{M}=\mathfrak{H}\left(\mathbf{x}_{\star}\right)\mathbf{H}_{\star}^{-1}$
we obtain
\begin{align}
\left[\left(\mu\mathbf{I}-\mathbf{W}\right)+\alpha\mathbf{W}\mathbf{M}\right]\left(\mu\mathbf{I}-\mathbf{W}\right)\mathbf{y} & =-\alpha\mathbf{W}\mathbf{M}\left(\mathbf{W}-\mathbf{I}\right)\mathbf{y} & \Rightarrow\nonumber \\
\left(\mu\mathbf{I}-\mathbf{W}\right)^{2}\mathbf{y}+\alpha\mathbf{W}\mathbf{M}\left(\mu\mathbf{I}-\mathbf{W}\right)\mathbf{y} & =\alpha\mathbf{W}\mathbf{M}\left(\mathbf{I}-\mathbf{W}\right)\mathbf{y} & \Rightarrow\nonumber \\
\left(\mu\mathbf{I}-\mathbf{W}\right)^{2}\mathbf{y} & =\alpha\left(1-\mu\right)\mathbf{W}\mathbf{M}\mathbf{y}.\label{eq:geigenv}
\end{align}
The above means that $\alpha\left(1-\mu\right)$ is a generalized
eigenvalue of the matrix pair $\left(\left(\mu\mathbf{I}-\mathbf{W}\right)^{2},\mathbf{W}\mathbf{M}\right)$.

Let $\mathbf{u}$ and $\mathbf{v}$ be, respectively, a right and
left eigenvector of $\mathbf{W}$, associated with $\mu_{0}$. We
have 
\[
\left(\mu_{0}\mathbf{I}-\mathbf{W}\right)^{2}\mathbf{u}=0.
\]

It then follows from Lemma~\ref{lem:geigenv} that 
\begin{align}
\alpha\left(1-\mu\right) & \simeq\frac{\mathbf{v}^{\top}\left(\mu\mathbf{I}-\mathbf{W}\right)^{2}\mathbf{u}}{\mathbf{v}^{\top}\mathbf{W}\mathbf{M}\mathbf{u}}\nonumber \\
 & =\frac{\mathbf{v}^{\top}\left[\left(\mu\mathbf{I}-\mathbf{W}\right)^{2}-\left(\mu_{0}\mathbf{I}-\mathbf{W}\right)^{2}\right]\mathbf{u}}{\mathbf{v}^{\top}\mathbf{W}\mathbf{M}\mathbf{u}}\nonumber \\
 & =\frac{2\left(\mu_{0}-\mu\right)\mathbf{v}^{\top}\mathbf{W}\mathbf{u}+\left(\mu^{2}-\mu_{0}^{2}\right)\mathbf{v}^{\top}\mathbf{u}}{\mathbf{v}^{\top}\mathbf{W}\mathbf{M}\mathbf{u}}.\label{eq:quad}
\end{align}
Clearly, $\mathbf{u}=u\otimes y$ and $\mathbf{v}=v\otimes z$, for
any $y,z\in\mathbb{R}^{N}$. Since we can choose $z$ arbitrarily,
we choose $z=y$. We then get
\begin{align}
\mathbf{v}^{\top}\mathbf{u} & =\left(v^{\top}u\right)\left(z^{\top}y\right)=v^{\top}u,\label{eq:aux1}\\
\mathbf{y}^{\top}\mathbf{W}\breve{\mathbf{x}} & =\mu_{0}\left(v^{\top}u\right)\left(z^{\top}y\right)=\mu_{0}v^{\top}u.\label{eq:aux2}
\end{align}
and
\begin{equation}
\mathbf{v}^{\top}\mathbf{W}\mathbf{M}\mathbf{u}=\mu_{0}y^{\top}\left(\sum_{i=1}^{I}v_{i}u_{i}\nabla^{2}f^{i}\left(x_{\star}\right)\left[\nabla^{2}f\left(x_{\star}\right)\right]^{-1}\right)y=\mu_{0}sv^{\top}u.\label{eq:aux3}
\end{equation}
The result follows by putting~(\ref{eq:aux1})-(\ref{eq:aux3}) into~(\ref{eq:quad})
and rearranging terms.
\end{IEEEproof}
We now use Theorem~\ref{thm:eigenapprox} to analyze the trajectory
of the relevant eigenvalues of $\Gamma(\alpha)$, when $\alpha$ is
small. 

The largest eigenvalue of $\boldsymbol{\Gamma}(0)$ is $1$. For this
case, we have $u=v=\mathbf{1}_{I}$. Then $R=\mathbf{I}_{N}$ and
$s=1$. Hence~(\ref{eq:sol}) becomes
\[
\mu^{2}-\left(2-\alpha\right)\mu+\left(1-\alpha\right)\simeq0,
\]
giving that either $\mu\simeq1$ or $\mu\simeq1-\alpha$. Since $\mu_{0}=1$
has multiplicity $2N$, this means that $\boldsymbol{\Gamma}(\alpha)$
will (approximately) have an eigenvalue at $\mu=1$, with multiplicity
$N$ and another one at $\mu=1-\alpha$ with the same multiplicity.
The first set of $N$ eigenvalues is consequence of the fact that
any point of the form $\left(\mathbf{x}_{k},\mathbf{g}_{k}\right)=\left(\mathbf{1}_{I}\otimes x,0\right)$,
for any $x\in\mathbb{R}^{N}$, is a stationary point of the local
linear dynamics determined by $\boldsymbol{\Gamma}(\alpha)$. However,
we know from the global nonlinear model~(\ref{eq:ssa1})-(\ref{eq:ssa3})
that the only possible of such stationary points is $\left(\mathbf{1}_{I}\otimes x_{\star},0\right)$.
Hence, the convergence speed is determined by the remaining $N(2I-1)$
eigenvalues of $\boldsymbol{\Gamma}(\alpha)$. Hence, the second set
of $N$ eigenvalues describes a convergence mode of the distributed
optimization algorithm

The second largest eigenvalue of $\boldsymbol{\Gamma}(0)$ is $\lambda_{2}$.
In this case we have
\[
\mu^{2}-\lambda_{2}\left(2-\alpha s\right)\mu+\lambda_{2}\left(\lambda_{2}-\alpha s\right)\simeq0
\]
which gives
\begin{align*}
\mu & \simeq\frac{\lambda_{2}}{2}\left(2-\alpha s\pm\sqrt{\alpha^{2}s^{2}+4\alpha s\left(\frac{1}{\lambda_{2}}-1\right)}\right).
\end{align*}
As before, the above means that $\boldsymbol{\Gamma}(\alpha)$ will
(approximately) have two eigenvalues with multiplicity $N$, one moving
up from $\lambda_{2}$ and another one moving down. Considering the
one moving up we can devise a criterion for choosing the design value
$\alpha_{\star}$ of $\alpha$. More precisely, we require this eigenvalue
to be equal to $1-\alpha_{\star}$, i.e.,
\begin{equation}
1-\alpha_{\star}=\left|\frac{\lambda_{2}}{2}\left(2-\alpha_{\star}s+\sqrt{\alpha_{\star}^{2}s^{2}+4\alpha_{\star}s\left(\frac{1}{\lambda_{2}}-1\right)}\right)\right|.\label{eq:opt-alpha}
\end{equation}

We now do the following approximation in~(\ref{eq:R}) 
\begin{equation}
R\simeq\frac{1}{I}\sum_{i=1}^{I}H^{i}\left(x_{\star}\right)H^{-1}\left(x_{\star}\right)=\mathbf{I}.\label{eq:approx}
\end{equation}
This gives $s=1$, which when put in equation~(\ref{eq:opt-alpha})
gives the following step size selection criterion
\begin{equation}
1-\alpha_{\star}=\left|\frac{\lambda_{2}}{2}\left(2-\alpha_{\star}+\sqrt{\alpha_{\star}^{2}+4\alpha_{\star}\left(\frac{1}{\lambda_{2}}-1\right)}\right)\right|.\label{eq:offline-alpha}
\end{equation}
In particular, if $\lambda_{2}\in\mathbb{R}$,
\[
\alpha_{\star}=1-\sqrt{\lambda_{2}}.
\]

\subsection{Distributed step size estimation\protect\label{subsec:Adaptive}}

The criterion~(\ref{eq:offline-alpha}) is based on the somehow coarse
approximation~(\ref{eq:approx}). A more accurate step size selection
can be achieved if estimates of $R$, $u$, $v$ and $\lambda_{2}$
are available at each node. The parameters $u$, $v$ and $\lambda_{2}$
depend on the communication network. These parameters can be either
known in advance, or estimated during an initialization stage. In
particular, if the network is undirected, the distributed method described
in Appendix~\ref{sec:Distributed-eig} can be used.

In contrast, matrix $R$ depends on the optimization problem. Hence
it needs to be estimated. We can do so in a distributed manner using
dynamic average consensus. Let $R_{k}^{i}$ denote the estimate of
$R$ obtained at node~$i$ and time $k$, and let $\mathbf{r}_{k}=\mathrm{col}\left(R_{k}^{1},\cdots,R_{k}^{I}\right)$.
The estimation then is initialized by $R_{1}^{i}=\mathbf{I}_{N}$,
for all $i=1,\cdots,I$, i.e., $\mathbf{r}_{1}=\mathrm{col}\left(\mathbf{I}_{N},\cdots,\mathbf{I}_{N}\right)$,
and proceeds as follows

\begin{equation}
\mathbf{r}_{k+1}=\mathbf{W}\left[\mathbf{r}_{k}+\mathfrak{r}\left(\mathbf{x}_{k+1}\right)-\mathfrak{r}\left(\mathbf{x}_{k}\right)\right],\label{eq:ssR}
\end{equation}
where $\mathfrak{r}\left(\mathbf{x}_{k}\right)=\mathrm{col}\left(\mathfrak{r}^{1}\left(\mathbf{x}_{k}\right),\cdots,\mathfrak{r}^{I}\left(\mathbf{x}_{k}\right)\right)$
with
\[
r^{i}\left(\mathbf{x}_{k}\right)=\frac{v^{i}u^{i}}{v^{\top}u}\nabla^{2}f^{i}\left(x_{k}^{i}\right)\left(H_{k}^{i}\right)^{-1}.
\]

In order to compute $\alpha_{k}^{i}$ at each node and time $k$ we
need to solve~(\ref{eq:opt-alpha}). This requires computing an approximation
$s_{k}^{i}$ of $s$ using $R_{k}^{i}$ in place of $R$. We do not
know the value of $y$ which gives the best approximation in~(\ref{eq:sol}).
But we know from~(\ref{eq:approx}) that $R_{k}^{i}$ approaches
$\mathbf{I}$ as $k$ increases. We then compute $s_{k}^{i}$ as the
mid point between the largest and smallest eigenvalues of $R_{k}^{i}$,
i.e., we choose 
\begin{equation}
s_{k}^{i}\simeq\frac{\left\Vert R_{k}^{i}\right\Vert +\left\Vert \left(R_{k}^{i}\right)^{-1}\right\Vert ^{-1}}{2}.\label{eq:adaptive-s}
\end{equation}
We then compute $\alpha_{k}^{i}$ by solving
\begin{equation}
1-\alpha_{k}^{i}=\left|\frac{\lambda_{2}}{2}\left(2-\alpha_{k}^{i}s_{k}^{i}+\sqrt{\left(\alpha_{k}^{i}s_{k}^{i}\right)^{2}+4\alpha_{k}^{i}s_{k}^{i}\left(\frac{1}{\lambda_{2}}-1\right)}\right)\right|.\label{eq:adaptive-alpha}
\end{equation}

The following theorem states the local stability of the distributed
optimization algorithm when used together with the adaptive step size
selection method.
\begin{thm}
\label{thm:global-1} Suppose Assumption~\ref{assu:local} holds
and, for all $k\in\mathbb{N}$ and $i=1,\cdots,I$, $\alpha_{k}^{i}$
is chosen using~(\ref{eq:adaptive-alpha}). Then, there exists a
neighborhood of $\left(\mathbf{x}_{\star},0,\mathbf{h}_{\star},\mathbf{r}_{\star}\right)$,
with $\mathbf{r}_{\star}=\mathbf{1}_{I}\otimes R$, such that if $\left(\mathbf{x}_{1},\mathbf{g}_{1},\mathbf{h}_{1},\mathbf{r}_{1}\right)$
is inside that neighborhood,
\begin{equation}
\lim_{k\rightarrow\infty}x_{k}^{i}=x_{\star},\quad\text{for all }i\in\{1,\dots,I\}.\label{eq:stability-1}
\end{equation}
\end{thm}
\begin{IEEEproof}
Let $\tilde{\mathbf{r}}_{k}=\mathbf{r}_{k}-\mathbf{r}_{\star}$. In
the adaptive step size selection algorithm described in Section~\ref{subsec:Adaptive},
$\boldsymbol{\alpha}_{k}$ is independent of $\boldsymbol{\xi}_{k}$
and only depends on $\mathbf{r}_{k}$. With some abuse of notation
we write $\boldsymbol{\alpha}_{k}=\boldsymbol{\alpha}\left(\mathbf{r}_{k}\right)$
and 
\begin{align}
\boldsymbol{\zeta}_{k+1} & =\mathfrak{F}\left(\boldsymbol{\xi}_{k},\boldsymbol{\alpha}\left(\tilde{\mathbf{r}}_{k}+\mathbf{r}_{\star}\right)\right),\label{eq:WSS1}\\
\tilde{\mathbf{r}}_{k+1} & =\tilde{\mathbf{W}}\left[\tilde{\mathbf{r}}_{k}+\mathfrak{r}\left(\mathbf{x}_{k+1}\right)-\mathfrak{r}\left(\mathbf{x}_{k}\right)\right],\label{eq:WSS2}
\end{align}
where $\mathfrak{F}$ represents the mapping induced by~(\ref{eq:ssa1})-(\ref{eq:ssa3}).
Clearly, $\left(\boldsymbol{\xi}_{k},\tilde{\mathbf{r}}_{k}\right)=\left(0,0\right)$
is an equilibrium point of the above system.

Now, since $\text{\ensuremath{\boldsymbol{\alpha}_{k}}}$ always appears
multiplying $\ensuremath{\mathbf{g}_{k}}$ in~(\ref{eq:ssa1})-(\ref{eq:ssa3}),
it straightforwardly follows that 
\[
\mathscr{D}_{\boldsymbol{\alpha}}\mathfrak{F}\left(0,\boldsymbol{\alpha}_{\star}\right)\left(\boldsymbol{\alpha}_{k}-\boldsymbol{\alpha}_{\star}\right)=0,
\]
where $\boldsymbol{\alpha}_{\star}=\boldsymbol{\alpha}\left(\mathbf{r}_{\star}\right)$.
Hence, in the local linear dynamics of~(\ref{eq:WSS1})-(\ref{eq:WSS2})
around $\left(0,\boldsymbol{\alpha}_{\star}\right)$, the term $\left(\boldsymbol{\alpha}_{k}-\boldsymbol{\alpha}_{\star}\right)$
does not act as input of $\boldsymbol{\xi}_{k+1}$. We also have,
from the discussion in Section~\ref{subsec:Adaptive}, that with
the step size choice $\boldsymbol{\alpha}_{\star}$,~(\ref{eq:WSS1})
is locally stable in a neighborhood of $\boldsymbol{\xi}_{k}=0$.
This means that $\mathfrak{F}\left(\boldsymbol{\xi}_{k},\boldsymbol{\alpha}_{\star}\right)$
is locally stable in that neighborhood. The local stability of~(\ref{eq:WSS1})-(\ref{eq:WSS2})
in the same neighborhood then immediately follows from the stability
of~(\ref{eq:WSS2}).
\end{IEEEproof}

\section{Step size selection for guaranteed global convergence\protect\label{sec:Step-size-global}}

In this section we consider the case in which the objective function
$f$ has a single local minimizer $x_{\star}$. We provide a time-varying
step size selection scheme to guarantee that, under certain regularity
conditions, the convergence~(\ref{eq:stability}) of the proposed
method.

In this section we do the following assumption.
\begin{assumption}
\label{assu:global} There exist constants $\gamma,\delta>0$ such
that 
\begin{align*}
\sup_{x\in\mathbb{R}^{N}}\left\Vert \left[\nabla^{2}f(x)\right]^{-1}\right\Vert  & \leq\beta, & \sup_{\substack{x\in\mathbb{R}^{N}\\
1\leq i\leq I
}
}\left\Vert \nabla^{2}f^{i}(x)\right\Vert  & \leq\gamma, & \sup_{\substack{x\in\mathbb{R}^{N}\\
1\leq i\leq I
}
}\frac{\left\Vert \nabla^{2}f^{i}(x)-\nabla^{2}f^{i}(y)\right\Vert _{\mathrm{F}}}{\left\Vert x-y\right\Vert } & \leq\delta.
\end{align*}
\end{assumption}
\begin{rem}
\label{rem:convex} Notice that the first condition of Assumption~\ref{assu:global}
requires that $f$ is strongly convex, which in turn implies the existence
of a single local minimum $x_{\star}$. This condition is a weaker
requirement than that of most global stability results from the literature
which, as mentioned in Section~\ref{sec:Introduction}, require strong
convexity of each local function $f^{i}$. Notice also that the Lipschitz
continuity assumption on the Hessian implies the existence of a finite
$\gamma$ when $x$ is restricted to any bounded subset of $\mathbb{R}^{N}$.
\end{rem}
Assumption~\ref{assu:global} implies the following properties.
\begin{lem}
\label{lem:handy}Under Assumption~\ref{assu:global}, 
\begin{align*}
\sup_{\mathbf{x}\neq\mathbf{y}}\frac{\left\Vert \mathfrak{g}(\mathbf{x})-\mathfrak{g}(\mathbf{y})\right\Vert }{\left\Vert \mathbf{x}-\mathbf{y}\right\Vert } & \leq\gamma, & \sup_{\mathbf{x}\neq\mathbf{y}}\frac{\left\Vert \mathfrak{H}(\mathbf{x})-\mathfrak{H}(\mathbf{y})\right\Vert _{\mathrm{F}}}{\left\Vert \mathbf{x}-\mathbf{y}\right\Vert } & \leq\delta, & \sup_{\mathbf{x}\neq\mathbf{y}}\frac{\left\Vert \bar{\mathfrak{H}}(\mathbf{x})-\bar{\mathfrak{H}}(\mathbf{y})\right\Vert _{\mathrm{F}}}{\left\Vert \mathbf{x}-\mathbf{y}\right\Vert } & \leq\delta.
\end{align*}
\end{lem}
\begin{IEEEproof}
We have
\begin{align*}
\left\Vert \mathfrak{g}\left(\mathbf{x}\right)-\mathfrak{g}\left(\mathbf{y}\right)\right\Vert ^{2} & =\sum_{i=1}^{I}\left\Vert \nabla f^{i}\left(x^{i}\right)-\nabla f^{i}\left(y^{i}\right)\right\Vert ^{2}\\
 & \leq\gamma^{2}\sum_{i=1}^{I}\left\Vert x^{i}-y^{i}\right\Vert ^{2}\\
 & =\gamma^{2}\left\Vert \mathbf{x}-\mathbf{y}\right\Vert ^{2}.
\end{align*}
Also
\begin{align*}
\left\Vert \mathfrak{H}\left(\mathbf{x}\right)-\mathfrak{H}\left(\mathbf{y}\right)\right\Vert _{\mathrm{F}}^{2} & =\sum_{i=1}^{I}\left\Vert \nabla^{2}f^{i}\left(x^{i}\right)-\nabla^{2}f^{i}\left(y^{i}\right)\right\Vert _{\mathrm{F}}^{2}\\
 & \leq\delta^{2}\sum_{i=1}^{I}\left\Vert x^{i}-y^{i}\right\Vert ^{2}\\
 & \leq\delta^{2}\left\Vert \mathbf{x}-\mathbf{y}\right\Vert ^{2},
\end{align*}
and
\begin{align*}
\left\Vert \bar{\mathfrak{H}}\left(\mathbf{x}\right)-\bar{\mathfrak{H}}\left(\mathbf{y}\right)\right\Vert _{\mathrm{F}} & =\sqrt{I}\left\Vert \frac{1}{I}\sum_{i=1}^{I}\nabla^{2}f^{i}\left(x^{i}\right)-\nabla^{2}f^{i}\left(y^{i}\right)\right\Vert _{\mathrm{F}}\\
 & \leq\frac{1}{\sqrt{I}}\sum_{i=1}^{I}\left\Vert \nabla^{2}f^{i}\left(x^{i}\right)-\nabla^{2}f^{i}\left(y^{i}\right)\right\Vert _{\mathrm{F}}\\
 & \leq\frac{\delta}{\sqrt{I}}\sum_{i=1}^{I}\left\Vert x^{i}-y^{i}\right\Vert \\
 & \leq\delta\left\Vert \mathbf{x}-\mathbf{y}\right\Vert .
\end{align*}
\end{IEEEproof}
We now introduce the notation required for stating our main result.
\begin{notation}
\label{not:bounds} Let $W=T^{-1}\Lambda T$ be the Jordan decomposition
of $W$ and $\mathbf{T}=T\otimes\mathbf{I}_{N}$. Let $\theta_{k}^{\top}=\left[\left\Vert \mathbf{T}\tilde{\mathbf{x}}_{k}\right\Vert ,\left\Vert \mathbf{T}\tilde{\mathbf{g}}_{k}\right\Vert ,\left\Vert \mathbf{T}\tilde{\mathbf{h}}_{k}\right\Vert _{\mathrm{F}}\right]$.
We define $\tau=\left\Vert T\right\Vert $, $\varsigma=\left\Vert T^{-1}\right\Vert $
as well as
\begin{align*}
\mu^{\top} & =\left[\varsigma\gamma,\varsigma,0\right],\\
\nu^{\top} & =\left[\varsigma\beta\delta\left(\beta\gamma+1\right),0,\varsigma\beta\right],\\
\psi^{\top} & =\left[\beta\tau,\beta\gamma\tau\upsilon,\beta\delta\tau\upsilon\right]
\end{align*}
and
\[
\Omega=\frac{\varsigma^{2}}{2}\left[\begin{array}{ccc}
\beta\gamma\delta\left(\beta\gamma+2\right) & \beta\delta & \beta\gamma\\
\beta\delta & \beta^{2}\delta & \beta\\
\beta\gamma & \beta & 0
\end{array}\right],\qquad\Phi=\left[\begin{array}{ccc}
\lambda_{2} & 0 & 0\\
\gamma\tau\eta\upsilon & \lambda_{2} & 0\\
\delta\tau\eta\upsilon & 0 & \lambda_{2}
\end{array}\right],\qquad\Psi=\left[\begin{array}{ccc}
\beta\gamma\tau & \beta\tau & 0\\
\beta\gamma^{2}\tau\upsilon & \beta\gamma\tau\upsilon & 0\\
\beta\gamma\delta\tau\upsilon & \beta\delta\tau\upsilon & 0
\end{array}\right],
\]
Let $P\in\mathbb{R}^{3\times3}$ be the unique positive solution of
$\Phi^{\top}P\Phi=P-I.$ Let $\eta=\left\Vert \mathbf{I}-W\right\Vert $,
$\upsilon=\left\Vert A-W\right\Vert $,
\begin{align*}
\mathsf{a} & =\mu^{\top}P^{-1}\mu, & \mathsf{b} & =\frac{\beta^{2}\delta}{2}, & \mathsf{c} & =\nu^{\top}P^{-1}\nu,\\
\mathsf{d} & =\left\Vert \Omega^{1/2}P^{-1/2}\right\Vert ^{2}, & \mathsf{e} & =\left\Vert P^{-1}\right\Vert , & \mathsf{h} & =\left\Vert \psi^{\top}P\psi\right\Vert .
\end{align*}
and
\begin{align*}
\mathsf{f}(\alpha) & =\left\Vert \Psi^{\top}P\left(2\Phi+\alpha\Psi\right)P^{-1}\right\Vert ,\\
\mathsf{g}(\alpha) & =2\left\Vert \psi^{\top}P\left(\Phi+k\Psi\right)P^{-1}\left(\Phi+\alpha\Psi\right)^{\top}P\psi\right\Vert .
\end{align*}
Finally, we define the mapping $F:\mathbb{R}^{2}\times\mathbb{R}\rightarrow\mathbb{R}^{2}:\left(\left(x_{1},x_{2}\right),\alpha\right)\mapsto\left(y_{1},y_{2}\right)$
by
\begin{align*}
y_{1} & =\left(1-\alpha\right)x_{1}+\alpha\mathsf{a}x_{2}+\alpha^{2}\mathsf{b}x_{1}^{2}+\alpha\mathsf{c}x_{1}x_{2}+\alpha\mathsf{d}x_{2}^{2},\\
y_{2}^{2} & =\left(1-\mathsf{e}+\alpha\mathsf{f}\left(\alpha_{k}\right)\right)x_{2}^{2}+\alpha\mathsf{g}\left(\alpha\right)x_{2}x_{1}+\alpha^{2}\mathsf{h}\left(k\right)x_{1}^{2}.
\end{align*}
\end{notation}
The following result gives the required global convergence condition.
\begin{thm}
\label{thm:global} Let $\alpha_{k}^{i}=\alpha_{k}$, for all $k\in\mathbb{N}$
and $i=1,\cdots,I$, where 
\begin{align*}
\alpha_{k} & =\underset{\alpha}{\arg\min}\left\Vert F\left(\chi_{k},\alpha\right)\right\Vert ,\\
\text{s.t.} & \chi_{k}\leq F\left(\chi_{k},\alpha\right).
\end{align*}
and $\chi_{k}$ generated by the following iterations 
\[
\chi_{k+1}=F\left(\chi_{k},\alpha_{k}\right),
\]
initialized by some $\chi_{1}^{\top}\geq\left[\left\Vert \bar{\mathfrak{g}}\left(\bar{\mathbf{x}}_{1}\right)\right\Vert ,\left\Vert \theta_{1}\right\Vert _{P}^{2}\right]$.
Then, under Assumption~\ref{assu:global}, 
\begin{equation}
\lim_{k\rightarrow\infty}x_{k}^{i}=x_{\star},\quad\text{for all }i\in\{1,\dots,I\}.\label{eq:stability}
\end{equation}
\end{thm}
\begin{rem}
Notice that the constraint $\chi_{k}\leq F\left(\chi_{k},\alpha\right)$
implies that
\[
\alpha_{k}\mathsf{a}\left\Vert \bar{\mathfrak{g}}\left(\bar{\mathbf{x}}_{k}\right)\right\Vert ^{2}+\mathsf{b}\left\Vert \theta_{k}\right\Vert _{P}+\mathsf{c}\left\Vert \bar{\mathfrak{g}}\left(\bar{\mathbf{x}}_{k}\right)\right\Vert \left\Vert \theta_{k}\right\Vert _{P}+\mathsf{d}\left\Vert \theta_{k}\right\Vert _{P}^{2}\leq\left\Vert \bar{\mathfrak{g}}\left(\bar{\mathbf{x}}_{k}\right)\right\Vert ,
\]
and that $\theta_{k}$ measures the inter node-node variable mismatch.
If $\theta_{k}$ is too large, it may occur that the above can only
be satisfied if $\alpha_{k}=0$. In such case, the algorithm automatically
runs a number of pure consensus iterations, i.e., chooses $\alpha_{k}=0$,
until the above inequality can be satisfied with $\alpha_{k}>0$.
\end{rem}
The rest of the section is devoted to show the above result.
\begin{lem}
\label{lem:GS1} If Assumption~\ref{assu:global} holds and $0\leq\alpha_{k}^{i}=\alpha_{k}\leq1$,
for all $i=1,\cdots,I$, then 
\[
\left\Vert \bar{\mathfrak{g}}\left(\bar{\mathbf{x}}_{k+1}\right)\right\Vert \leq\left(1-\alpha_{k}\right)\left\Vert \bar{\mathfrak{g}}\left(\bar{\mathbf{x}}_{k}\right)\right\Vert +\alpha_{k}\mathsf{a}\left\Vert \theta_{k}\right\Vert _{P}+\alpha_{k}^{2}\mathsf{b}\left\Vert \bar{\mathfrak{g}}\left(\bar{\mathbf{x}}_{k}\right)\right\Vert ^{2}+\alpha_{k}\mathsf{c}\left\Vert \bar{\mathfrak{g}}\left(\bar{\mathbf{x}}_{k}\right)\right\Vert \left\Vert \theta_{k}\right\Vert _{P}+\alpha_{k}\mathsf{d}\left\Vert \theta_{k}\right\Vert _{P}^{2},
\]
\end{lem}
\begin{IEEEproof}
Using the following inequality, 
\begin{align}
\left\Vert \mathbf{g}_{k}\right\Vert  & \leq\left\Vert \bar{\mathfrak{g}}\left(\bar{\mathbf{x}}_{k}+\tilde{\mathbf{x}}_{k}\right)\right\Vert +\left\Vert \tilde{\mathbf{g}}_{k}\right\Vert \nonumber \\
 & \leq\left\Vert \bar{\mathfrak{g}}\left(\bar{\mathbf{x}}_{k}\right)\right\Vert +\left\Vert \bar{\mathfrak{g}}\left(\bar{\mathbf{x}}_{k}+\tilde{\mathbf{x}}_{k}\right)-\bar{\mathfrak{g}}\left(\bar{\mathbf{x}}_{k}\right)\right\Vert +\left\Vert \tilde{\mathbf{g}}_{k}\right\Vert \nonumber \\
 & \leq\left\Vert \bar{\mathfrak{g}}\left(\bar{\mathbf{x}}_{k}\right)\right\Vert +\gamma\left\Vert \tilde{\mathbf{x}}_{k}\right\Vert +\left\Vert \tilde{\mathbf{g}}_{k}\right\Vert ,\label{eq:g-bound}
\end{align}
as well as the equality $\mathbf{A}\mathfrak{H}\left(\bar{\mathbf{x}}_{k}\right)\mathbf{A}\mathbf{y}=\bar{\mathfrak{H}}\left(\bar{\mathbf{x}}_{k}\right)\mathbf{A}\mathbf{y}$,
for any $\mathbf{y}\in\mathbb{R}^{IN}$, we obtain
\begin{align}
\bar{\mathfrak{g}}\left(\bar{\mathbf{x}}_{k+1}\right) & =\bar{\mathfrak{g}}\left(\bar{\mathbf{x}}_{k}-\alpha_{k}\mathbf{A}\mathbf{B}_{k}^{-1}\mathbf{g}_{k}\right)\nonumber \\
 & =\bar{\mathfrak{g}}\left(\bar{\mathbf{x}}_{k}\right)-\left[\int_{0}^{1}\mathbf{A}\mathfrak{H}\left(\bar{\mathbf{x}}_{k}-t\alpha_{k}\mathbf{A}\mathbf{B}_{k}^{-1}\mathbf{g}_{k}\right)dt\right]\alpha_{k}\mathbf{A}\mathbf{B}_{k}^{-1}\mathbf{g}_{k}\nonumber \\
 & =\bar{\mathfrak{g}}\left(\bar{\mathbf{x}}_{k}\right)-\alpha_{k}\bar{\mathfrak{H}}\left(\bar{\mathbf{x}}_{k}\right)\mathbf{A}\mathbf{B}_{k}^{-1}\mathbf{g}_{k}+\epsilon_{k}^{(1)}\nonumber \\
 & =\bar{\mathfrak{g}}\left(\bar{\mathbf{x}}_{k}\right)-\alpha_{k}\bar{\mathfrak{g}}\left(\bar{\mathbf{x}}_{k}\right)+\alpha_{k}\mathbf{A}\bar{\mathfrak{g}}\left(\bar{\mathbf{x}}_{k}\right)-\alpha_{k}\mathbf{A}\mathbf{g}_{k}+\alpha_{k}\mathbf{A}\mathbf{B}_{k}\mathbf{B}_{k}^{-1}\mathbf{g}_{k}-\alpha_{k}\mathbf{A}\bar{\mathfrak{H}}\left(\bar{\mathbf{x}}_{k}\right)\mathbf{B}_{k}^{-1}\mathbf{g}_{k}+\epsilon_{k}^{(1)}\nonumber \\
 & =\left(1-\alpha_{k}\right)\bar{\mathfrak{g}}\left(\bar{\mathbf{x}}_{k}\right)+\epsilon_{k}^{(1)}+\epsilon_{k}^{(2)}+\epsilon_{k}^{(3)},\label{eq:g-bar}
\end{align}
with
\begin{align*}
\epsilon_{k}^{(1)} & =\left\{ \int_{0}^{1}\left[\bar{\mathfrak{H}}\left(\bar{\mathbf{x}}_{k}\right)-\bar{\mathfrak{H}}\left(\bar{\mathbf{x}}_{k}-t\mathbf{A}\alpha_{k}\mathbf{B}_{k}^{-1}\mathbf{g}_{k}\right)\right]dt\right\} \mathbf{A}\alpha_{k}\mathbf{B}_{k}^{-1}\mathbf{g}_{k},\\
\epsilon_{k}^{(2)} & =\alpha_{k}\mathbf{A}\left(\bar{\mathfrak{g}}\left(\bar{\mathbf{x}}_{k}\right)-\mathbf{g}_{k}\right),\\
\epsilon_{k}^{(3)} & =\alpha_{k}\mathbf{A}\left(\mathbf{B}_{k}-\bar{\mathfrak{H}}\left(\bar{\mathbf{x}}_{k}\right)\right)\mathbf{B}_{k}^{-1}\mathbf{g}_{k}.
\end{align*}

Now
\begin{align}
\left\Vert \epsilon_{k}^{(1)}\right\Vert  & \leq\int_{0}^{1}\left\Vert \bar{\mathfrak{H}}\left(\bar{\mathbf{x}}_{k}-t\alpha_{k}\mathbf{A}\mathbf{B}_{k}^{-1}\mathbf{g}_{k}\right)-\bar{\mathfrak{H}}\left(\bar{\mathbf{x}}_{k}\right)\right\Vert dt\left\Vert \alpha_{k}\mathbf{B}_{k}^{-1}\mathbf{g}_{k}\right\Vert \nonumber \\
 & \leq\frac{\beta^{2}\delta\alpha_{k}^{2}}{2}\left\Vert \mathbf{g}_{k}\right\Vert ^{2}\nonumber \\
 & \leq\frac{\beta^{2}\delta\alpha_{k}^{2}}{2}\left(\left\Vert \bar{\mathfrak{g}}\left(\bar{\mathbf{x}}_{k}+\tilde{\mathbf{x}}_{k}\right)\right\Vert ^{2}+\left\Vert \tilde{\mathbf{g}}_{k}\right\Vert ^{2}\right)\nonumber \\
 & \leq\frac{\beta^{2}\delta\alpha_{k}^{2}}{2}\left(\left\Vert \bar{\mathfrak{g}}\left(\bar{\mathbf{x}}_{k}\right)\right\Vert ^{2}+\gamma^{2}\left\Vert \tilde{\mathbf{x}}_{k}\right\Vert ^{2}+2\gamma\left\Vert \bar{\mathfrak{g}}\left(\bar{\mathbf{x}}_{k}\right)\right\Vert \left\Vert \tilde{\mathbf{x}}_{k}\right\Vert +\left\Vert \tilde{\mathbf{g}}_{k}\right\Vert ^{2}\right).\label{eq:ep1}
\end{align}
Also
\begin{align}
\left\Vert \epsilon_{k}^{(2)}\right\Vert  & \leq\alpha_{k}\left\Vert \bar{\mathfrak{g}}\left(\bar{\mathbf{x}}_{k}\right)-\mathbf{g}_{k}\right\Vert \nonumber \\
 & \leq\alpha_{k}\left(\left\Vert \tilde{\mathbf{g}}_{k}\right\Vert +\left\Vert \bar{\mathfrak{g}}\left(\bar{\mathbf{x}}_{k}+\tilde{\mathbf{x}}_{k}\right)-\bar{\mathfrak{g}}\left(\bar{\mathbf{x}}_{k}\right)\right\Vert \right)\nonumber \\
 & \leq\alpha_{k}\left(\left\Vert \tilde{\mathbf{g}}_{k}\right\Vert +\gamma\left\Vert \tilde{\mathbf{x}}_{k}\right\Vert \right),\label{eq:ep2}
\end{align}

Since $\bar{\mathfrak{H}}\left(\bar{\mathbf{x}}_{k}\right)=\mathrm{diag}\left(\nabla^{2}f(\bar{x}),\cdots,\nabla^{2}f(\bar{x})\right)$,
we have $\left\Vert \mathfrak{B}\left(\bar{\mathfrak{H}}\left(\bar{\mathbf{x}}_{k}\right)\right)-\bar{\mathfrak{H}}\left(\bar{\mathbf{x}}_{k}\right)\right\Vert _{\mathrm{F}}=0$.
We also have from~\cite{wihler2009holder} that $\left\Vert \mathfrak{B}(\mathbf{x})-\mathfrak{B}(\mathbf{y})\right\Vert _{\mathrm{F}}\leq\left\Vert \mathbf{x}-\mathbf{y}\right\Vert $,
for all $\mathbf{x},\mathbf{y}\in\mathbb{R}^{IN}$. We then obtain
\begin{align*}
\left\Vert \mathbf{B}_{k}-\bar{\mathfrak{H}}\left(\bar{\mathbf{x}}_{k}\right)\right\Vert  & \leq\left\Vert \mathfrak{B}\left(\bar{\mathfrak{H}}\left(\mathbf{x}_{k}\right)+\tilde{\mathbf{H}}_{k}\right)-\mathfrak{B}\left(\bar{\mathfrak{H}}\left(\mathbf{x}_{k}\right)\right)\right\Vert _{\mathrm{F}}\\
 & +\left\Vert \mathfrak{B}\left(\bar{\mathfrak{H}}\left(\mathbf{x}_{k}\right)\right)-\mathfrak{B}\left(\bar{\mathfrak{H}}\left(\bar{\mathbf{x}}_{k}\right)\right)\right\Vert _{\mathrm{F}}\\
 & +\left\Vert \mathfrak{B}\left(\bar{\mathfrak{H}}\left(\bar{\mathbf{x}}_{k}\right)\right)-\bar{\mathfrak{H}}\left(\bar{\mathbf{x}}_{k}\right)\right\Vert _{\mathrm{F}}\\
 & \leq\left\Vert \tilde{\mathbf{H}}_{k}\right\Vert _{\mathrm{F}}+\left\Vert \bar{\mathfrak{H}}\left(\mathbf{x}_{k}\right)-\bar{\mathfrak{H}}\left(\bar{\mathbf{x}}_{k}\right)\right\Vert _{\mathrm{F}}\\
 & \leq\left\Vert \tilde{\mathbf{h}}_{k}\right\Vert _{\mathrm{F}}+\delta\left\Vert \tilde{\mathbf{x}}_{k}\right\Vert .
\end{align*}
Hence
\begin{align}
\left\Vert \epsilon_{k}^{(3)}\right\Vert  & \leq\beta\alpha_{k}\left\Vert \mathbf{B}_{k}-\bar{\mathfrak{H}}\left(\bar{\mathbf{x}}_{k}\right)\right\Vert \left\Vert \mathbf{g}_{k}\right\Vert \nonumber \\
 & \leq\beta\alpha_{k}\left(\left\Vert \tilde{\mathbf{H}}_{k}\right\Vert +\left\Vert \bar{\mathfrak{H}}\left(\bar{\mathbf{x}}_{k}+\tilde{\mathbf{x}}_{k}\right)-\bar{\mathfrak{H}}\left(\bar{\mathbf{x}}_{k}\right)\right\Vert \right)\left(\left\Vert \bar{\mathfrak{g}}\left(\bar{\mathbf{x}}_{k}\right)\right\Vert +\gamma\left\Vert \tilde{\mathbf{x}}_{k}\right\Vert +\left\Vert \tilde{\mathbf{g}}_{k}\right\Vert \right)\nonumber \\
 & \leq\beta\alpha_{k}\left(\left\Vert \tilde{\mathbf{h}}_{k}\right\Vert _{\mathrm{F}}+\delta\left\Vert \tilde{\mathbf{x}}_{k}\right\Vert \right)\left(\left\Vert \bar{\mathfrak{g}}\left(\bar{\mathbf{x}}_{k}\right)\right\Vert +\gamma\left\Vert \tilde{\mathbf{x}}_{k}\right\Vert +\left\Vert \tilde{\mathbf{g}}_{k}\right\Vert \right).\label{eq:ep3}
\end{align}
Taking the norm in~(\ref{eq:g-bar}), and using~(\ref{eq:ep1})-(\ref{eq:ep3}),
we obtain 
\[
\left\Vert \bar{\mathfrak{g}}\left(\bar{\mathbf{x}}_{k+1}\right)\right\Vert \leq\left(1-\alpha_{k}\right)\left\Vert \bar{\mathfrak{g}}\left(\bar{\mathbf{x}}_{k}\right)\right\Vert +\alpha_{k}^{2}\mathsf{b}\left\Vert \bar{\mathfrak{g}}\left(\bar{\mathbf{x}}_{k}\right)\right\Vert ^{2}+\alpha_{k}\left(\mu^{\top}+\nu^{\top}\left\Vert \bar{\mathfrak{g}}\left(\bar{\mathbf{x}}_{k}\right)\right\Vert +\theta_{k}^{\top}\Omega\right)\theta_{k}
\]
The result then follows since $\left\Vert \boldsymbol{y}\right\Vert \leq\varsigma\left\Vert \mathbf{T}\mathbf{y}\right\Vert $,
for $\mathbf{y}=\tilde{\mathbf{x}}_{k},\tilde{\mathbf{g}}_{k},\tilde{\mathbf{h}}_{k}$
and
\begin{align*}
y^{\top}x & =y^{\top}P^{-1}y\left\Vert x\right\Vert _{P},\\
\theta_{k}^{\top}\Omega\theta_{k} & =\left\Vert \Omega^{1/2}P^{-1/2}\right\Vert ^{2}\left\Vert x\right\Vert _{P}^{2}.
\end{align*}
\end{IEEEproof}
\begin{lem}
\label{lem:GS2} If Assumption~\ref{assu:global} holds and $0\leq\alpha_{k}^{i}=\alpha_{k}$,
for all $i=1,\cdots,I$, then
\[
\left\Vert \theta_{k+1}\right\Vert _{P}^{2}\leq\left\Vert \theta_{k}\right\Vert _{P}^{2}-\mathsf{e}\left\Vert \theta_{k}\right\Vert ^{2}+\alpha_{k}\mathsf{f}\left(\alpha_{k}\right)\left\Vert \theta_{k}\right\Vert ^{2}+\alpha_{k}\mathsf{g}\left(\alpha_{k}\right)\left\Vert \theta_{k}\right\Vert \left\Vert \bar{\mathfrak{g}}\left(\bar{\mathbf{x}}_{k}\right)\right\Vert +\alpha_{k}^{2}\mathsf{h}\left(\alpha_{k}\right)\left\Vert \bar{\mathfrak{g}}\left(\bar{\mathbf{x}}_{k}\right)\right\Vert ^{2}.
\]
\end{lem}
\begin{IEEEproof}
We start by bounding each entry of the vector $\theta_{k+1}$. Using
Lemma~\ref{lem:handy} we obtain
\begin{align*}
\left\Vert \mathbf{T}\tilde{\mathbf{g}}_{k+1}\right\Vert  & \leq\lambda_{2}\left\Vert \mathbf{T}\tilde{\mathbf{g}}_{k}\right\Vert +\tau\upsilon\left\Vert \tilde{\mathfrak{g}}\left(\bar{\mathbf{x}}_{k}+\mathbf{W}\tilde{\mathbf{x}}_{k}-\alpha_{k}\mathbf{B}_{k}^{-1}\mathbf{g}_{k}\right)-\tilde{\mathfrak{g}}\left(\bar{\mathbf{x}}_{k}+\tilde{\mathbf{x}}_{k}\right)\right\Vert \\
 & \leq\lambda_{2}\left\Vert \tilde{\mathbf{g}}_{k}\right\Vert +\gamma\tau\upsilon\left\Vert \left(\mathbf{I}-\mathbf{W}\right)\tilde{\mathbf{x}}_{k}-\alpha_{k}\mathbf{B}_{k}^{-1}\mathbf{g}_{k}\right\Vert \\
 & \leq\lambda_{2}\left\Vert \tilde{\mathbf{g}}_{k}\right\Vert +\gamma\tau\upsilon\eta\left\Vert \tilde{\mathbf{x}}_{k}\right\Vert +\alpha_{k}\gamma\beta\tau\upsilon\left\Vert \mathbf{g}_{k}\right\Vert \\
 & \leq\lambda_{2}\left\Vert \tilde{\mathbf{g}}_{k}\right\Vert +\gamma\tau\upsilon\eta\left\Vert \tilde{\mathbf{x}}_{k}\right\Vert +\alpha_{k}\gamma\beta\tau\upsilon\left(\left\Vert \bar{\mathfrak{g}}\left(\bar{\mathbf{x}}_{k}\right)\right\Vert +\gamma\left\Vert \tilde{\mathbf{x}}_{k}\right\Vert +\left\Vert \tilde{\mathbf{g}}_{k}\right\Vert \right)\\
 & \leq\alpha_{k}\beta\gamma\tau\upsilon\left\Vert \bar{\mathfrak{g}}\left(\bar{\mathbf{x}}_{k}\right)\right\Vert +\gamma\tau\upsilon\left(\eta+\alpha_{k}\beta\gamma\right)\left\Vert \tilde{\mathbf{x}}_{k}\right\Vert +\left(\lambda_{2}+\alpha_{k}\beta\gamma\tau\upsilon\right)\left\Vert \tilde{\mathbf{g}}_{k}\right\Vert 
\end{align*}
Following similar steps we obtain
\begin{align*}
\left\Vert \mathbf{T}\tilde{\mathbf{h}}_{k+1}\right\Vert _{\mathrm{F}} & \leq\alpha_{k}\beta\delta\tau\upsilon\left\Vert \bar{\mathfrak{g}}\left(\bar{\mathbf{x}}_{k}\right)\right\Vert +\delta\tau\upsilon\left(\eta+\alpha_{k}\beta\gamma\right)\left\Vert \tilde{\mathbf{x}}_{k}\right\Vert +\alpha_{k}\beta\delta\tau\upsilon\left\Vert \tilde{\mathbf{g}}_{k}\right\Vert +\lambda_{2}\left\Vert \tilde{\mathbf{h}}_{k}\right\Vert _{\mathrm{F}}
\end{align*}
and
\begin{align*}
\left\Vert \mathbf{T}\tilde{\mathbf{x}}_{k+1}\right\Vert  & \leq\lambda_{2}\left\Vert \mathbf{T}\tilde{\mathbf{x}}_{k}\right\Vert +\alpha_{k}\tau\left\Vert \mathbf{B}_{k}^{-1}\right\Vert \left\Vert \mathbf{g}_{k}\right\Vert \\
 & \leq\alpha_{k}\beta\tau\left\Vert \bar{\mathfrak{g}}\left(\bar{\mathbf{x}}_{k}\right)\right\Vert +\left(\lambda_{2}+\alpha_{k}\beta\gamma\tau\right)\left\Vert \tilde{\mathbf{x}}_{k}\right\Vert +\alpha_{k}\beta\tau\left\Vert \tilde{\mathbf{g}}_{k}\right\Vert .
\end{align*}

From the above we obtain 
\[
\theta_{k+1}=\left(\Phi+\alpha_{k}\Psi\right)\theta_{k}+\alpha_{k}\psi\left\Vert \bar{\mathfrak{g}}\left(\bar{\mathbf{x}}_{k}\right)\right\Vert .
\]
Then
\begin{align*}
\left\Vert \theta_{k+1}\right\Vert _{P}^{2} & =\theta_{k}^{\top}\left(\Phi+\alpha_{k}\Psi\right)^{\top}P\left(\Phi+\alpha_{k}\Psi\right)\theta_{k}+2\alpha_{k}\psi^{\top}P\left(\Phi+\alpha_{k}\Psi\right)\theta_{k}\left\Vert \bar{\mathfrak{g}}\left(\bar{\mathbf{x}}_{k}\right)\right\Vert +\alpha_{k}^{2}\psi^{\top}P\psi\left\Vert \bar{\mathfrak{g}}\left(\bar{\mathbf{x}}_{k}\right)\right\Vert ^{2}\\
 & =\theta_{k}^{\top}\Phi^{\top}P\Phi\theta_{k}+2\alpha_{k}\theta_{k}^{\top}\Phi^{\top}P\Psi\theta_{k}+\alpha_{k}^{2}\theta_{k}^{\top}\Psi^{\top}P\Psi\theta_{k}\\
 & +\alpha_{k}\psi^{\top}P\left(2\left(\Phi+\alpha_{k}\Psi\right)\theta_{k}+\alpha_{k}\psi\left\Vert \bar{\mathfrak{g}}\left(\bar{\mathbf{x}}_{k}\right)\right\Vert \right)\left\Vert \bar{\mathfrak{g}}\left(\bar{\mathbf{x}}_{k}\right)\right\Vert \\
 & =\theta_{k}^{\top}P\theta_{k}-\theta_{k}^{\top}\theta_{k}+\alpha_{k}\theta_{k}^{\top}\left(2\Phi^{\top}+\alpha_{k}\Psi^{\top}\right)P\Psi\theta_{k}\\
 & +\alpha_{k}\psi^{\top}P\left(2\left(\Phi+\alpha_{k}\Psi\right)\theta_{k}+\alpha_{k}\psi\left\Vert \bar{\mathfrak{g}}\left(\bar{\mathbf{x}}_{k}\right)\right\Vert \right)\left\Vert \bar{\mathfrak{g}}\left(\bar{\mathbf{x}}_{k}\right)\right\Vert \\
 & \leq\left\Vert \theta_{k}\right\Vert _{P}^{2}-\left\Vert P^{-1}\right\Vert \left\Vert \theta_{k}\right\Vert _{P}^{2}+\alpha_{k}\left\Vert \Psi^{\top}P\left(2\Phi+\alpha_{k}\Psi\right)P^{-1}\right\Vert \left\Vert \theta_{k}\right\Vert _{P}^{2}\\
 & +2\alpha_{k}\left\Vert \psi^{\top}P\left(\Phi+\alpha_{k}\Psi\right)P^{-1}\left(\Phi+\alpha_{k}\Psi\right)^{\top}P\psi\right\Vert \left\Vert \theta_{k}\right\Vert _{P}\left\Vert \bar{\mathfrak{g}}\left(\bar{\mathbf{x}}_{k}\right)\right\Vert +\alpha_{k}^{2}\left\Vert \psi^{\top}P\psi\right\Vert \left\Vert \bar{\mathfrak{g}}\left(\bar{\mathbf{x}}_{k}\right)\right\Vert ^{2},
\end{align*}
and the result follows.
\end{IEEEproof}
\begin{IEEEproof}[Proof of Theorem~\ref{thm:global}]
 Let $\zeta_{k}^{\top}=\left[\left\Vert \bar{\mathfrak{g}}\left(\bar{\mathbf{x}}_{k}\right)\right\Vert ,\left\Vert \theta_{k}\right\Vert _{P}\right]$.
It follows from Lemmas~\ref{lem:GS1} and~\ref{lem:GS2} that, for
any $0\leq\alpha\leq1$, 
\[
\zeta_{k+1}\leq F\left(\zeta_{k},\alpha\right).
\]
The proof proceed by induction. We have $\zeta_{1}\leq\chi_{1}$.
At iteration $k$, suppose $\zeta_{k}\leq\chi_{k}$. Since the map
$\chi\mapsto F(\chi,\alpha)$ is monotonous, we have $\zeta_{k+1}\leq\chi_{k+1}$.
The result would then follow if $\lim_{k\rightarrow\infty}\chi_{k}=0$.
This in turn occurs if both components of $\chi_{k}$ are strictly
monotonously decreasing. Notice that the constraint $F\left(\chi_{k},\alpha\right)\leq\chi_{k}$
is equivalent to 
\begin{align*}
\alpha_{k}\mathsf{a}\left\Vert \bar{\mathfrak{g}}\left(\bar{\mathbf{x}}_{k}\right)\right\Vert ^{2}+\mathsf{b}\left\Vert \theta_{k}\right\Vert _{P}+\mathsf{c}\left\Vert \bar{\mathfrak{g}}\left(\bar{\mathbf{x}}_{k}\right)\right\Vert \left\Vert \theta_{k}\right\Vert _{P}+\mathsf{d}\left\Vert \theta_{k}\right\Vert _{P}^{2} & \leq\left\Vert \bar{\mathfrak{g}}\left(\bar{\mathbf{x}}_{k}\right)\right\Vert ,\\
\alpha_{k}\mathsf{f}\left(\alpha_{k}\right)\left\Vert \theta_{k}\right\Vert _{P}^{2}+\alpha_{k}\mathsf{g}\left(\alpha_{k}\right)\left\Vert \theta_{k}\right\Vert _{P}\left\Vert \bar{\mathfrak{g}}\left(\bar{\mathbf{x}}_{k}\right)\right\Vert +\alpha_{k}^{2}\mathsf{h}\left(\alpha_{k}\right)\left\Vert \bar{\mathfrak{g}}\left(\bar{\mathbf{x}}_{k}\right)\right\Vert ^{2} & \leq\mathsf{e}\left\Vert \theta_{k}\right\Vert _{P}^{2}.
\end{align*}
Hence, there always exists $\alpha$ such that $\left[\chi_{k+1}\right]_{2}<\left[\chi_{k}\right]_{2}$.
Also, for $\left[\chi_{k+1}\right]_{1}<\left[\chi_{k}\right]_{1}$
to hold we need
\[
\left(\mathsf{b}+\mathsf{c}\left\Vert \bar{\mathfrak{g}}\left(\bar{\mathbf{x}}_{k}\right)\right\Vert +\mathsf{d}\left\Vert \theta_{k}\right\Vert _{P}\right)\left\Vert \theta_{k}\right\Vert _{P}\leq\left\Vert \bar{\mathfrak{g}}\left(\bar{\mathbf{x}}_{k}\right)\right\Vert .
\]
The convergence to zero of the first component is then guaranteed
by that of the second one.
\end{IEEEproof}

\section{A numerical example\protect\label{sec:Numerical-example}}

\subsection{Case study}

We consider a target localization problem. There are $I=30$ nodes,
measuring the distance to a target located at $x_{\mathrm{true}}\in\mathbb{R}^{2}$.
Node $i$ is located at $a^{i}\in\mathbb{R}^{2}$, with $a^{i}\sim\mathcal{N}\left(x_{\mathrm{true}},100\times\mathbf{I}_{2}\right)$,
and is initialized by $x_{1}^{i}\sim\mathcal{N}\left(x_{\mathrm{true}},\mathbf{I}_{2}\right)$.
It measures 
\[
z^{i}=\left\Vert x_{\mathrm{true}}-a^{i}\right\Vert ^{2}+n^{i},\qquad n^{i}\sim\mathcal{N}\left(0,\sigma^{2}\right),
\]
with $\sigma^{2}=0.01$. Nodes are connected via a network with ring
topology, whose gains are given by
\[
w^{i,j}=\begin{cases}
0.7, & i=j,\\
0.15, & \mod\left(i-j,I\right)=1,\\
0.15, & \mod\left(i-j,I\right)=I-2,\\
0, & \text{otherwise},
\end{cases}
\]
where $\mod\left(a,b\right)$ denotes the $a$ modulo $b$ operation.
This results in $\lambda_{2}=0.9838$. 

Doing a maximum likelihood estimation of $x$ we obtain
\begin{align*}
x_{\star} & =\underset{x}{\arg\max}p\left(z^{i},\cdots,z^{I}|x\right)\\
 & =\underset{x}{\arg\min}\sum_{i=1}^{I}f^{i}(x),
\end{align*}
with
\[
f^{i}(x)=\left(\left\Vert x-a^{i}\right\Vert ^{2}-z^{i}\right)^{2}.
\]
It is straightforward to obtain
\begin{align*}
\nabla f^{i}(x) & =4\left(\left\Vert x-a^{i}\right\Vert ^{2}-z^{i}\right)\left(x-a^{i}\right),\\
\nabla^{2}f^{i}(x) & =8\left(x-a^{i}\right)\left(x-a^{i}\right)^{\top}\\
 & +4\left(\left\Vert x-a^{i}\right\Vert ^{2}-z^{i}\right)\mathbf{I}_{2},
\end{align*}
Using the above, the optimal step size, in the sense of minimizing
the largest modulus of the eigenvalues of $\boldsymbol{\Gamma}(\alpha)$,
which differs from $1$, is given by $\alpha_{\mathrm{opt}}=6.197\times10^{-3}$.

\subsection{Numerical experiments}

\begin{figure}[h]
\begin{centering}
\begin{minipage}[t]{0.49\columnwidth}%
\begin{center}
\includegraphics[width=1\columnwidth]{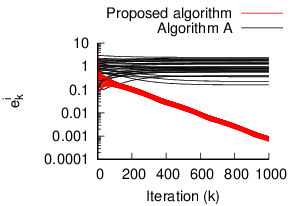}
\par\end{center}
\begin{center}
(a)
\par\end{center}%
\end{minipage}%
\begin{minipage}[t]{0.49\columnwidth}%
\begin{center}
\includegraphics[width=1\columnwidth]{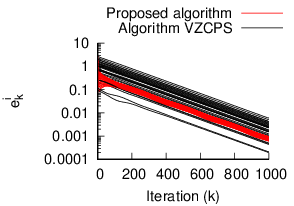}
\par\end{center}
\begin{center}
(b)
\par\end{center}%
\end{minipage}
\par\end{centering}
\caption{Effect of modification~(\ref{eq:algA}), i.e, removing consensus
on variables. (a) Algorithm~A does not converge. (b) Algorithm~VZCPS,
which also introduces modification~(\ref{eq:algB1})-(\ref{eq:algB2})
converges slower than the proposed one.}
\label{fig:close}
\end{figure}

In the first experiment we evaluate the effect of considering the
modification~(\ref{eq:algA}) introduced by Algorithms~A, as described
in Section~\ref{sec:Comparison-with-similar}. We use $x_{\mathrm{true}}=[0,0]^{\top}$,
$\beta=0.1$ and the optimal step size $\alpha_{\mathrm{opt}}=6.197\times10^{-3}$.
We also use $e_{k}^{i}\triangleq\left\Vert x_{k}^{i}-x_{\star}\right\Vert $
as the performance metric for each node. In Figure~\ref{fig:close}-(a)
we compare the performance of the proposed algorithm with Algorithm~A.
We see how the lack of a consensus stage prevents the local variables
$x_{k}^{i}$ to converge to a common value. In Figure~\ref{fig:close}-(b)
we see the effect of considering also the modification~(\ref{eq:algB1})-(\ref{eq:algB2}),
i.e, Algorithm~VZCPS. We see that it converges, although at a much
smaller rate than the proposed algorithm. We then conclude that it
is this second modification which causes the converge of Algorithm~VZCPS.

\begin{figure}[h]
\begin{centering}
\begin{minipage}[t]{0.33\columnwidth}%
\begin{center}
\includegraphics[width=1\columnwidth]{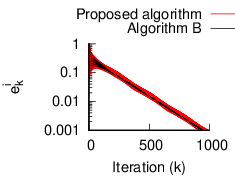}
\par\end{center}
\begin{center}
(a)
\par\end{center}%
\end{minipage}%
\begin{minipage}[t]{0.33\columnwidth}%
\begin{center}
\includegraphics[width=1\columnwidth]{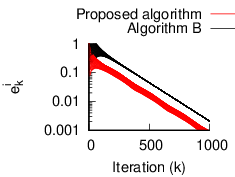}
\par\end{center}
\begin{center}
(b)
\par\end{center}%
\end{minipage}%
\begin{minipage}[t]{0.33\columnwidth}%
\begin{center}
\includegraphics[width=1\columnwidth]{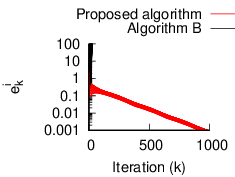}
\par\end{center}
\begin{center}
(c)
\par\end{center}%
\end{minipage}
\par\end{centering}
\caption{Effect of modification~(\ref{eq:algB1})-(\ref{eq:algB2}). (a) With
$x_{\mathrm{true}}=[0,0]^{\top}$ Algorithm~B performs similar to
the proposed one. (b) With $x_{\mathrm{true}}=[300,300]^{\top}$ the
parameters are pulled away from $x_{\star}$ until consensus is reached.
(c) With $x_{\mathrm{true}}=[1000,1000]^{\top}$ the pulling effect
is intensified causing instability.}
\label{fig:far}
\end{figure}

In the second experiment we remove modification~(\ref{eq:algA}),
i.e, add consensus on variables, and study the effect of modification~(\ref{eq:algB1})-(\ref{eq:algB2})
introduced by Algorithm~B. In Figure~\ref{fig:far}-(a) we see that
Algorithm~B converges at rate similar to that of the proposed algorithm.
However, as explained in Section~\ref{sec:proposed-method}, modification~(\ref{eq:algB1})-(\ref{eq:algB2})
has a negative effect when the minimizing parameters $x_{\star}$
are far from zero. We show this in Figure~\ref{fig:far}-(b), where
we repeat the previous experiment with $x_{\mathrm{true}}=[300,300]^{\top}$.
We see how the local estimates of Algorithm~B are pulled away from
$x_{\star}$ during the initial iterations, until consensus is reached.
In Figure~\ref{fig:far}-(c) we repeat the experiment with $x_{\mathrm{true}}=[1000,1000]^{\top}$.
We see that Algorithm~B is not able to reach consensus in time, which
causes its divergence.

\begin{figure}[h]
\begin{centering}
\begin{minipage}[t]{0.33\columnwidth}%
\begin{center}
\includegraphics[width=1\columnwidth]{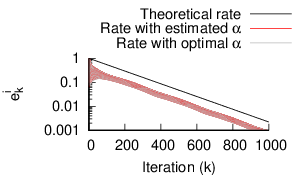}
\par\end{center}
\begin{center}
(a)
\par\end{center}%
\end{minipage}%
\begin{minipage}[t]{0.33\columnwidth}%
\begin{center}
\includegraphics[width=1\columnwidth]{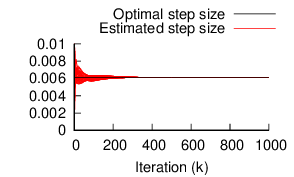}
\par\end{center}
\begin{center}
(b)
\par\end{center}%
\end{minipage}%
\begin{minipage}[t]{0.33\columnwidth}%
\begin{center}
\includegraphics[width=1\columnwidth]{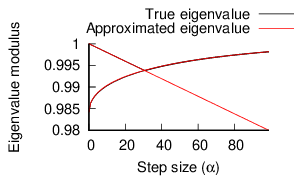}
\par\end{center}
\begin{center}
(c)
\par\end{center}%
\end{minipage}
\par\end{centering}
\caption{Distributed estimation of the optimal step size $\alpha_{\star}$.
(a) The asymptotic convergence rate matches the theoretically optimal
one. (b) Evolution of the distributed step size estimates. (c) Comparison
between the actual eigenvalue dependence on the step size and the
approximated one.}
\label{fig:adaptive}
\end{figure}

In the third experiment we evaluate the use of the distributed algorithm
for estimating the step size. In Figure~\ref{fig:adaptive}-(a) we
compare the convergence of $\left\Vert x_{k}^{i}-x_{\star}\right\Vert $
using both, the optimal step size $\alpha_{\mathrm{opt}}$ and the
distributedly estimated one. We see that both methods converge at
a very similar rate. We also show in the figure the theoretically
optimal rate $\rho^{k}\left(\boldsymbol{\Phi}-\alpha_{\mathrm{opt}}\boldsymbol{\Psi}\right)$.
We see that the asymptotic convergence rate of both methods closely
resembles the theoretical one. We also show in Figure~\ref{fig:adaptive}-(b)
the evolution of the estimated step size $\alpha_{k}^{i}$ at each
node. Finally, Figure~\ref{fig:adaptive}-(c) shows how the two eigenvalues
used to compute $\alpha_{\star}$ depend on $\alpha$, and compares
this with the approximated dependence given by Theorem~\ref{thm:eigenapprox}.
We see how, before the two eigenvalues meet, the true and approximated
trajectories closely resemble each other. This results in $\alpha_{\star}=6.117\times10^{-3}$
being a good approximation of $\alpha_{\mathrm{opt}}=6.197\times10^{-3}$.

\section{Conclusion\protect\label{sec:Conclusion}}

We aimed at achieving the fastest convergence rate for distributed
optimization. We did two steps towards this goal. In the first step
we proposed a new distributed optimization method which converges
faster than other available options. In the second step we proposed
a distributed method to estimate the step size that maximizes this
rate. We provided sufficient conditions for the convergence of the
resulting method in a neighborhood of a local solution. We also provided
condition to guarantee global convergence of the method, in the case
of an objective function having a single local minimum, together with
a different step size selection strategy. We present numerical experiments
confirming our claims.

\appendices{}

\section{Fréchet Derivatives\protect\label{subsec:Fr=0000E9chet-Derivatives}}

The Fréchet derivative generalizes the concept of derivative of functions
between Euclidean spaces to functions between normed vector spaces~\cite{Hamilton1982,Coleman2012}.
\begin{defn}
Let $X$ and $Z$ be normed vector spaces and $U\subseteq X$ be open.
A function $f:U\rightarrow Z$ is called Fréchet differentiable at
$x\in U$ if there exists a bounded linear map $A:X\rightarrow Z$
such that 
\[
\lim_{\left\Vert h\right\Vert \rightarrow0}\frac{\left\Vert f(x+h)-f(x)-A(h)\right\Vert }{\left\Vert h\right\Vert }=0.
\]
In this case we say that $A$ is the Fréchet derivative of $f$ at
$x$, and denote it by $\mathscr{D}f(x)=A$. We also say that $f$
is Fréchet differentiable if it is so at all $x\in U$. We use $\mathscr{D}f:x\mapsto\mathscr{D}f(x)$
to denote the Fréchet derivative of $f$ and $\mathscr{D}:f\mapsto\mathscr{D}f$
to denote the Fréchet derivation operator.

The $n$-th order Fréchet derivative is the $n$-fold composition
of the Fréchet derivation operator, i.e., 
\[
\mathscr{D}^{n}f=\underbrace{\mathscr{D}\circ\cdots\circ\mathscr{D}}_{t\text{ times}}f.
\]

If $f:X\times Y\rightarrow Z$, we define the partial Fréchet derivative
$\mathscr{D}_{x}f(x,y)$ of $f$ with respect to $x$ at $(x,y)$,
as the Fréchet derivative of the map $x\mapsto f(x,y)$ at $x$.
\end{defn}
As with the derivative of functions between Euclidean spaces, we can
approximate a function between normed spaces using a Taylor expansion.
\begin{lem}[Taylor theorem for Fréchet derivatives]
\label{lem:taylor-frechet} If $f:X\rightarrow Z$ is $n+1$-times
continuously differentiable, then 
\begin{align*}
f(x+h) & =f(x)+\mathscr{D}f(x)(h)+\frac{1}{2!}\mathscr{D}^{2}f(x)(h)(h)+\\
 & \cdots+\frac{1}{n!}\mathscr{D}^{n}f(x)(h)\cdots(h)+O\left(\left\Vert h\right\Vert ^{n+1}\right).
\end{align*}
\end{lem}

\section{Generalized eigenvalue problem}

Given a matrix pair $A,B\in\mathbb{R}^{N\times N}$, the generalized
eigenvalue problem consists in finding the values of $\lambda\in\mathbb{C}$,
satisfying 
\[
\det\left(A-\lambda B\right)=0.
\]
We say that $x,y\in\mathbb{R}^{N}$ are, respectively, right and left
generalized eigenvectors associated with $\lambda$ if 
\begin{align*}
Ax & =\lambda Bx,\\
y^{\top}A & =\lambda y^{\top}B.
\end{align*}

The following result gives an approximation of the perturbation $\tilde{\lambda}$
of a generalized eigenvalue $\lambda$, when matrix $A$ is modified
by adding to it a perturbation matrix $\tilde{A}$. 
\begin{lem}
\label{lem:geigenv}Let $x$ and $y$ be right and left generalized
eigenvectors of $A,B\in\mathbb{R}^{N\times N}$ associated with $\lambda$
and 
\[
\left(A+\tilde{A}\right)\left(x+\tilde{x}\right)=\left(\lambda+\tilde{\lambda}\right)B\left(x+\tilde{x}\right).
\]
Then
\[
\tilde{\lambda}=\frac{y^{\top}\tilde{A}\left(x+\tilde{x}\right)}{y^{\top}B\left(x+\tilde{x}\right)}.
\]
\end{lem}
\begin{IEEEproof}
We have
\begin{align*}
\left(A+\tilde{A}\right)\left(x+\tilde{x}\right) & =\left(\lambda+\tilde{\lambda}\right)B\left(x+\tilde{x}\right) & \Rightarrow\\
Ax+\tilde{A}x+A\tilde{x}+\tilde{A}\tilde{x} & =\lambda Bx+\tilde{\lambda}Bx+\lambda B\tilde{x}+\tilde{\lambda}B\tilde{x} & \Rightarrow\\
\tilde{A}x+A\tilde{x}+\tilde{A}\tilde{x} & =\tilde{\lambda}Bx+\lambda B\tilde{x}+\tilde{\lambda}B\tilde{x} & \Rightarrow\\
\tilde{A}\left(x+\tilde{x}\right)+A\tilde{x} & =\tilde{\lambda}B\left(x+\tilde{x}\right)+\lambda B\tilde{x}
\end{align*}
Then
\begin{align*}
0 & =y^{\top}\tilde{A}\left(x+\tilde{x}\right)+y^{\top}A\tilde{x}-\tilde{\lambda}y^{\top}B\left(x+\tilde{x}\right)-\lambda y^{\top}B\tilde{x}\\
 & =y^{\top}\tilde{A}\left(x+\tilde{x}\right)+y^{\top}\left(A-\lambda B\right)\tilde{x}-\tilde{\lambda}y^{\top}B\left(x+\tilde{x}\right)\\
 & =y^{\top}\tilde{A}\left(x+\tilde{x}\right)-\tilde{\lambda}y^{\top}B\left(x+\tilde{x}\right),
\end{align*}
and the result follows.
\end{IEEEproof}

\section{Distributed estimation $\lambda_{2}$ and $u$ in the case of undirected
communication graphs\protect\label{sec:Distributed-eig}}

In this section we describe a distributed method for estimating $\lambda_{2}$
and $u$ in the case where the graph induced by the communication
network is undirected.

We know that $1$ is the largest eigenvalue of $W$ with (left and
right) eigenvector $\mathbf{1}_{I}$. Hence, 
\[
V=W-\frac{1}{I}\mathbf{1}_{I}\mathbf{1}_{I}^{\top},
\]
has $\lambda_{2}$ as its largest eigenvalue with eigenvector $v$.
We can then obtain recursive estimates $\hat{\lambda}_{2,k}$ and
$\hat{u}_{k}$ , of $\lambda_{2}$ and $u$, using the power method.
To this end, an initialization random vector $\hat{u}_{1}$ is produced
by locally drawing each random entry at each node, and the following
iterations are run 
\begin{align*}
\hat{u}_{k+1} & =V\hat{u}_{k}\\
 & =W\hat{u}_{k}-\mathbf{1}_{I}\bar{u}.
\end{align*}
In order to run the above, we need a recursive estimate $\hat{\bar{u}}_{k}$
of $\bar{u}$ at each node. We also obtain so using the power method,
i.e., we put $\hat{\bar{u}}_{1}=\hat{u}_{1}$ and run 
\[
\hat{\bar{u}}_{k+1}=W\hat{\bar{u}}_{k}.
\]
The estimation of $\lambda_{2}$ is then obtained as follows 
\[
\hat{\lambda}_{2,k}=\frac{\left\Vert \hat{u}_{k}\right\Vert }{\left\Vert \hat{u}_{k-1}\right\Vert }.
\]

\bibliographystyle{unsrt}
\bibliography{references}

\end{document}